\documentclass[smallcondensed,numbook]{svjour3}
\smartqed
\journalname{Applied Categorical Structures}

\usepackage{enumerate,graphicx}
\usepackage[mathscr]{eucal}
\usepackage{amssymb,amscd}

\usepackage{amsmath}
\usepackage{latexsym}
\usepackage{mathptmx}
\usepackage{verbatim,url}

\input xymatrix
\input xyarrow
\usepackage{xy,diagrams11}
\usepackage{tikz}
\usetikzlibrary{matrix,arrows}

\newtheorem{thm}{Theorem}[section]
\newtheorem{lem}[thm]{Lemma}
\newtheorem{prop}[thm]{Proposition}
\newtheorem{cor}[thm]{Corollary}

\newcommand{\lincol}{black}
\newcommand{\linth}{thick}
\newcommand{\po}[2][\pocol]
{\filldraw[#1](#2) circle (2 pt);}
\newcommand{\li}[1]{\draw[\linth,\lincol] #1;}

\newenvironment{newlist}
   {\begin{list}{}{\setlength{\labelsep}{0.25cm}
                   \setlength{\labelwidth}{0.65cm}
                   \setlength{\leftmargin}{0.9cm}}}
   {\end{list}}

\newcommand{\defn}[1]{{\emph{#1}}}
\newcommand{\eusbA}{\medsub e {\kern-0.75pt\A\kern-0.75pt}}

\newcommand{\cat}[1]{\boldsymbol{\mathscr{#1}}}
\newcommand{\CA}{{\cat A}}
\newcommand{\CB}{{\cat B}}
\newcommand{\CCD}{\cat D}
\newcommand{\CX}{\cat X}
\newcommand{\CY}{\cat Y}
\newcommand{\CL}{\cat L}
\newcommand{\CG}{\cat G}

\font\bmi=cmmi8 scaled 1440
\newcommand{\powerset}{\raise.6ex\hbox{\bmi\char'175 }}
\font\intix=cmex9
\newcommand{\DM}{{\hbox{$\textfont3=\intix\displaystyle{\bigwedge}$}}}
\newcommand{\dJ}{{{\hbox{$\textfont3=\intix\displaystyle{\bigvee}$}}}}

\newcommand{\mph}[2]{\CL^{\rm mp}(#1,#2)}
\newcommand{\sph}[2]{\CL^{\rm sp}(#1,#2)} 
\newcommand{\mpe}[2]{\CG^{\rm mp}(#1,#2)}
\newcommand{\mpm}[2]{\CG_\T^{\rm mp}(#1,#2)}

\newcommand{\D}[1]{\mathrm D(#1)}

\newcommand{\Dflat}[1]{\mathrm{D}^\flat(#1)}
\newcommand{\fC}[1]{\mathrm C(#1)} 

\newcommand{\sub}[1]{_{_{\kern-.9pt{\scriptstyle #1}}}}
\newcommand{\medsub}[2]{#1\lower0.6ex\hbox{$\scriptstyle{#2}$}}
\newcommand{\eA}[1]{\medsub e {\kern-0.75pt\A\kern-0.75pt}(#1)}
\newcommand{\esub}[1]{\medsub e {\kern-0.75pt #1 \kern-0.75pt}}
\newcommand{\epsub}[1]{\medsub \varepsilon {\kern-1.25pt #1}}
\newcommand{\esubA}{\medsub e {\kern-0.75pt\A\kern-0.75pt}}

\newcommand{\bbar}[1]{{\underline{\mathbf{#1}}}}
\newcommand{\twiddle}[1]{{\smash{\underset{\raise.375ex\hbox{$\smash\sim$}}
       {\mathbf{#1}}}\vphantom{\underline{\mathbf{#1}}}}} 
\newcommand{\stwiddle}[1]{\smash{\underset{\smash{\raise.1ex\hbox{\small$\sim$}}}
                         {\mathbf{#1}}}\vphantom{#1}}
\newcommand{\MT}{\twiddle M}

\newcommand{\twoB}{\bbar 2}

\newcommand{\twoT}{\twiddle 2}

\newcommand{\T}{\mathscr{T}}

\newcommand{\A}{\mathbf A}

\newcommand{\Lalg}{\mathbf L}
\newcommand{\K}{\mathbf K}

\newcommand{\M}{\mathbf M}

\newcommand{\X}{\mathbf X}
\newcommand{\Y}{\mathbf Y}

\newcommand{\Ys}{\Y\kern -2pt _s}

\newcommand{\Tp}{{\mathscr{T}}}
\newcommand{\dom}{dom}
\newcommand{\ISP}{\mathbb{ISP}}

\newcommand{\bigand}{\mathop{\bigwedge\kern -8.5truept \bigwedge}}
\newcommand{\Bigand}{\mathop{\bigwedge\kern -10truept \bigwedge}}
\newcommand{\littleand}{\mathbin{\wedge\kern -8truept \wedge}}
\newcommand{\bigor}{\mathop{\bigvee\kern -8.5truept \bigvee}}
\newcommand{\Bigor}{\mathop{\bigvee\kern -10truept \bigvee}}
\newcommand{\littleor}{\mathbin{\vee\kern -8truept \vee}}
\renewcommand{\le}{\leqslant}
\renewcommand{\ge}{\geqslant}

\newcommand{\lsem}{[\kern-1.75pt[}
\newcommand{\rsem}{]\kern-1.75pt]}

\begin{document}

\title{A fresh perspective on canonical extensions for bounded lattices}

\titlerunning{Canonical extensions for bounded lattices}

\author{A.\,P.\,K. Craig         \and
        M. Haviar \and
        H.\,A. Priestley
				\thanks{The first author gratefully acknowledges funding from the Rhodes Trust. The second author acknowledges support from Slovak grants VEGA 1/0485/09
				and APVV-0223-10.}
}

\institute{A.\,P.\,K. Craig \at
              Mathematical Institute, University of Oxford, 24/29 St.~Giles,
              Oxford OX1 3LB, United Kingdom \\
              \email{craig@maths.ox.ac.uk}
           \and
           M. Haviar \at
              Faculty of Natural Sciences, Matej Bel University, Tajovsk\'eho 40,
              974 01 Bansk\'{a} Bystrica, Slovak Republic\\
              \email{miroslav.haviar@umb.sk}
           \and
           H.\,A. Priestley \at
           Mathematical Institute, University of Oxford, 24/29 St~Giles,
           Oxford OX1 3LB, United Kingdom\\
           \email{hap@maths.ox.ac.uk}
}

\date{Received: date / Accepted: date}

\maketitle

\begin{abstract}
This paper presents a novel treatment of the canonical extension of a bounded lattice, in the spirit of the
theory of natural dualities. At the level of objects, this can be achieved by exploiting  the topological
representation due to M.~Plo\v{s}\v{c}ica, and the canonical extension can be obtained in the same manner
as can be done in the distributive case by exploiting Priestley duality. To encompass both objects and morphisms
the  Plo\v{s}\v{c}ica representation is replaced by a duality due to Allwein and Hartonas, recast in the style of Plo\v{s}\v{c}ica's paper. This leads to a  construction of canonical extension valid for all bounded lattices,
which is shown to be functorial, with the property that the canonical extension functor decomposes as
the composite of two functors, each of which acts on morphisms by composition, in the manner of hom-functors.
\end{abstract}

\section{Introduction} \label{intro}

Two valuable tools have been developed to study lattice-based algebras:
the theory of canonical extensions and the theory of topological dualities.
Although the two methodologies  have largely been developed separately,
they have a symbiotic relationship. In this paper we present a new approach
to the canonical extension of a bounded lattice, based on the dual representation due to
Plo\v{s}\v{c}ica~\cite{Pl95}. This yields a construction with strong
affinities with the theory of natural dualities---one which has clear categorical merits.

We shall consider only lattices which are bounded, and include the universal
 bounds, $0$ and~$1$, in the signature; we denote the
variety of all such lattices by~$\CL$. The subvariety of $\CL$ consisting of distributive lattices is denoted by
$\CCD$. The canonical extension of a member $\Lalg$ of $\CL$ is a
particular lattice completion $\Lalg^\delta $ of $\Lalg$,
having  two  properties, known as
\defn{density} and \defn{compactness}; we recall the definitions
below. These properties characterise this completion  uniquely, up to
an isomorphism fixing the original lattice.

The existence of canonical extensions
for members of~$\CL$ was first established by Gehrke
and Harding \cite{GH01}, drawing on the theory of Galois connections;
specifically,~$\Lalg^\delta $ arises as the complete lattice of
Galois-stable
sets associated with the polarity $R$ between the
filter lattice $\text{Filt} (\Lalg)$ and the ideal lattice
$\text{Idl} (\Lalg)$ of $\Lalg$ given by
$(F,I)\in R$
if and only if $F \cap I \ne
\emptyset$.  Further insight into the relationship between $\Lalg$,
$\text{Filt}(\Lalg)$,  $\text{Idl}(\Lalg)$ and $\Lalg^\delta $ was
supplied by Gehrke and Priestley \cite{GP08}.

The original investigation  by Gehrke and J\'onsson \cite{GJ94}
of canonical extensions in the distributive case
built $\Lalg^\delta $ by exploiting Priestley duality. Of course,
uniqueness of the canonical extension ensures that
this construction leads to the same completion
as do the (specialisations of) the constructions in \cite{GH01} and \cite{GP08}.
However the  translations between the various concrete realisations
of the canonical extension
of a distributive lattice
are rather indirect.
More significantly, there are
obstacles to extending, in a
fully satisfactory and transparent
way, the  duality approach beyond the distributive case.
Therefore we contend that it is worthwhile to explore
in  greater detail than hitherto the interface between
canonical extensions and duality theory for lattices,
and a categorical context for these.

Each of Priestley duality for $\CCD$ and Stone duality for
Boolean algebras, $\CB$,
is an
instance of a natural duality. Indeed, these dualities
 provided prototypical examples for the general theory,
as presented in the text of Clark and Davey \cite{CD98}.
In its basic form,
the theory of natural dualities applies to
a prevariety $\CA = \ISP(\M)$, where $\M$ is some finite algebra.
One sets up  functors
$\mathrm{D} \colon \CA \to \CX_\T$ and $\mathrm{E}\colon \CX_\T \to \CA$
between $\CA$ and a suitable
category $\CX_\T$ of structured Boolean spaces.  Here $\mathrm{D}$ and $\mathrm{E}$ are
 hom-functors
into, respectively, $\M$ and
$\MT_\T$, where $\MT_\T$ is an alter ego for $\M$,
that is, an object in $\CX_\T$ whose
underlying set is the same as that of $\M$.
If, for  a suitable choice of $\MT_\T$, we obtain  $\A \cong \mathrm{ED}(\A)$
for all
$\A \in \CA$, then we say we have
a \defn{duality}; this suffices to give a topological representation for
the members of $\CA$.  (Optimally, one also has
$\X \cong \mathrm{DE}(\X)$
for all $\X \in \CX_\T$ and hence a
dual category equivalence between $\CA$ and $\CX_\T$, but such
\defn{full dualities}
are not of relevance here.)
In the case of $\CB$ and $\CCD$, the algebra $\M$ is a two-element
algebra, and $\MT_\T$ is the set $\{ 0,1\}$ equipped with the discrete topology
and,  respectively, no relations and a single relation, {\it viz.\/}
the order relation $\le$
with $0 < 1$. The functors~$\mathrm D$ and~$\mathrm E$ act on morphisms by
composition, and this is a key feature of the way in which the duality
operates.  Moreover, by modifying the functors ~$\mathrm D$ and~$\mathrm E$,
 one arrives at a functorial
construction of the canonical extension---one which makes overt
its categorical relationship to the duality.
The canonical extension construction can then be viewed in the
 way  shown in Fig.~\ref{fig:Dsquare}.  In the figure,
$\flat$ denotes the functor forgetting the topology and
the functor $\mathrm{G}$ acts in the same manner as~$\mathrm{E}$ does,
but on untopologised rather than topologised
structures. The codomain $\CCD^+$ of~$\mathrm{G}$ can be taken to be the category of
doubly algebraic
distributive lattices with complete lattice homomorphisms.
In fact, more is true:  there is a hom-functor~$\mathrm{F}$
adjoint to~$\mathrm{G}$
such that~$\mathrm{F}$ and~$\mathrm{G}$
set up a dual equivalence
 between~$\CCD^+$ and~$\CX$ (see~\cite{CE1,DHvP11}).
However the functor~$\mathrm{F}$ is not involved in the factorisation of~$\delta$.
On the other hand, we do make use of~$\mathrm{E}$:
\[
\Lalg \cong \mathrm{ED}(\Lalg) \subseteq \Lalg^\delta  \quad \text{for
all } \Lalg \in \CCD.
\]

\begin{figure}  [ht]
\begin{diagram}[labelstyle=\textstyle]
\CCD& \pile{\rTo^{\mathrm D}\\ \lTo_{\mathrm E}} & \CX_\Tp \\
\dTo^{^\delta} & & \dTo_{\flat}\\
\CCD^+&  \lTo_{\mathrm G} & \CX
\end{diagram}
\caption{Factorising the canonical extension functor
on~$\CCD$
\label{fig:Dsquare}}
\end{figure}

Unlike $\CB$ and $\CCD$, the
variety $\CL$ of bounded lattices
is not finitely generated.  This
prevents it
from  coming  within the scope of
the theory of natural dualities.
However there is available
a long-established
representation of bounded lattices due to Urquhart~\cite{U78}.
This does not stem  from a dual category equivalence
but it does provide a concrete
topological representation of each $\Lalg \in \CL$.
Of greater relevance to us here, though, is
the  recasting
of Urquhart's work  by Plo\v{s}\v{c}ica~\cite{Pl95}. This
was motivated by natural duality theory, and is in the
spirit of that theory.  In barest outline, Plo\v{s}\v{c}ica's extension of the Priestley
representation is accomplished by replacing total maps into $\{ 0,1\}$,
{\it qua}
lattice or {\it qua}   partially ordered set,
by appropriate maximally-defined partial maps of the same sort.  One
may then conjecture
that, {\it mutatis mutandis\/}, the canonical extension of an arbitrary
bounded lattice may be obtained in the same manner as in Fig.~\ref{fig:Dsquare}
by `forgetting the topology' at the level of the Plo\v{s}\v{c}ica first dual.
The principal result of Section~\ref{newCE},
Theorem~\ref{canext},  confirms that this is indeed the case.
In preparation, Section \ref{prelim}   summarises
the key  notions
presented in  Plo\v{s}\v{c}ica's paper \cite{Pl95},
 investigates
 the basic properties of the corresponding concepts for the
topology-free setting,
 and reveals the interplay
between the topologised and untopologised versions.
We emphasise that
the
validation of our
construction
of the canonical extension
is independent of pre-existing theory of canonical extensions for
$\CL$: we rely solely
on the characteristic density and compactness properties,
to ensure that, for $\Lalg \in \CL$,  our candidate for $\Lalg^\delta$
contains an isomorphic copy of~$\Lalg$.
(A detailed discussion
of the role of polarities in our framework,
and a reconciliation of alternative approaches will be presented
separately (Craig and Haviar,~\cite{CH}).)

So, at the object level, we have a construction which mimics
the duality technnique used in the distributive case.  But it
is with morphisms that obstacles
arise with duality theory for~$\CL$.   We are nonetheless able to present a functorial
construction, and one in which the canonical extension functor factorises
in a manner analogous to that shown in Fig.~\ref{fig:Dsquare}.
Some sacrifice  is unavoidable to achieve this, and in our approach
this is the enlargement of the first dual space of~$\Lalg$.
Specifically,
our factorisation  aligns with the duality for~$\CL$ due to Allwein and
Hartonas~\cite{A92,AH93} and so
it is expedient to recast the
Allwein and Hartonas representation at the object level in
the Plo\v{s}\v{c}ica style.
We retain a critical feature of the natural
 duality framework, {\it viz.\/} that the functors involved
behave very like hom-functors. By saying a functor, $\mathrm H$ say,
is `very like' a hom-functor, we mean
that there is some fixed target,~$\mathrm N$ say (not necessarily in the
codomain of~$\mathrm H$), such that~$\mathrm H$ sends an object~$C$ to a specified
 set of partial maps from~$C$
into~$\mathrm{N}$
and that~$\mathrm H$ acts on morphisms by  composition.
Thus we are able to retain key features present in the distributive
case at the cost of amending the ancillary functors into which
the functor $^\delta$ splits, but without having to change the very natural
way in which these functors act
(see Figure \ref{fig:Lsquare} in Section \ref{cat}).
The major work involved in carrying out our programme is in
verifying that the functors we require are well-defined.
Once this is done,
the lifting of homomorphisms to canonical extensions  goes
through very cleanly.
We could have worked with the Allwein and Hartonas duality from
the outset and essentially bypassed Plo\v s\v cica's representation.
We have elected not to do this for two reasons.
Firstly,  in certain applications morphisms are not important,
and it is then advantageous to have a construction which coincides
in the distributive case
with the usual duality-based one.
Secondly, we are able to
convert the Plo\v s\v cica-style construction at the object level to one
based on the enlarged first dual space, and so do not need to
start afresh to encompass morphisms.

We note that
canonical extensions of  lattices provide a platform
for studying lattice-based algebras (also known as lattice
expansions).  A valuable recent survey  of canonical extensions
for lattice-based algebras, including a discussion of their important role
in the semantic modelling of logics, is given by
Gehrke and Vosmaer~\cite{GV11}.

\section{The framework for the construction}  \label{prelim}

In this section we set up the framework within which we shall construct the canonical
extension of a bounded lattice $\Lalg$.  In outline, we look at the structures,
and associated maps, obtained by deleting the topology from the dual
spaces employed in  Plo\v s\v cica's representation. We thereby
arrive at a complete lattice which we shall later show serves as
$\Lalg^\delta $.

We begin by recalling some basic definitions concerning completions.
A \defn{completion} of a (bounded)
lattice $\Lalg$ is defined to be a pair $(e,C)$ where $C$ is a
complete lattice  and $e \colon L \hookrightarrow C$ is an embedding.
An element of a completion
$(e,C)$ of a (bounded) lattice~$L$ which is
representable as a meet (join) of elements from $e(L)$ is
called a \defn{filter element} (\defn{ideal element}).
Filter (ideal) elements are called closed (open) elements in the older literature.
The sets of such elements
will be  denoted
$\mathbb {F} (C)$ and $\mathbb {I} (C)$, respectively.
A completion $(e,C)$ of~$L$
is said to be \defn{dense} if every element of~$C$ is both a join
of meets and a meet of joins of elements from $e(L)$; it is said
to be \defn{compact} if, for any sets
$A,B \subseteq L$
with $\bigwedge e(A) \le \bigvee e(B)$,
there exist finite subsets $A' \subseteq A$ and $B' \subseteq B$
such that
$\bigwedge e(A') \le \bigvee e(B')$.
A \defn{canonical extension} of $\Lalg$ is a dense and
compact completion.
As noted above,
every bounded lattice $\Lalg$ has a
canonical extension and any two canonical extensions of $\Lalg$ are
isomorphic via an isomorphism that fixes the elements of~$\Lalg$.

The central idea in Plo\v s\-\v ci\-ca's representation of bounded lattices~\cite{Pl95}
is the replacement of total maps by partial maps.
Let $\Lalg_1, \Lalg_2 \in \CL$.
A partial map $f \colon \Lalg_1\to \Lalg_2$ between bounded lattices
is called a \defn{partial homomorphism} if its domain is a
0,1-sublattice of $\Lalg_1$ and
$f\colon \dom(f) \to \Lalg_2$
is an $\CL$-homomorphism.
A partial homomorphism is said to be
\defn{maximal}
if there is no partial homomorphism
properly extending it; such a map is referred to as an MPH for short.
By Zorn's Lemma, every partial homomorphism can be
extended to an MPH. For
bounded lattices $\Lalg$ and $\K$, we
denote by
$\mph{\Lalg}{\K}$ the set of all MPH's from $\Lalg$ to $\K$.

Let
\[
\twoB :=\langle \{0,1\}; \vee, \wedge, 0, 1\rangle \quad \text{and} \quad
 \twoT
:=\langle \{0,1\}; \le\rangle
 \]
denote, respectively, the two-element bounded lattice and the
two-element ordered set with $0 < 1$.
The topological structures $\twoB_\T$ and
$\twoT_\T$
are obtained by adding the discrete topology $\T$ to $\twoB$
and~$\twoT$, respectively.

Following Plo\v s\v cica \cite{Pl95},
for any bounded lattice $\Lalg$, the topological dual space of~$\Lalg$
is defined
in the following way.
We equip
the set $\mph{\Lalg}{\twoB}$ with
the binary relation $E$ defined by the rule
\[
(f,g)\in E \quad \text{iff}\ f(x)\le g(x)\ \text{for every }
x\in\dom(f)\cap\dom(g);
\]
when needed, we also endow it with the topology $\T$ which has a
subbasis of closed sets all sets of the form
\[
V_a=\{\,f \in \mph{\Lalg}{\twoB} \mid f(a)=0 \,\} \quad \text{and} \quad
W_a=\{\,f \in \mph{\Lalg}{\twoB} \mid f(a)=1 \,\},
\]
where $a\in L$.
We let $\D{\Lalg} = (\mph{\Lalg}{\twoB}, E, \T)$. The topology of $\D{\Lalg}$ is T$_1$ and moreover,
it is compact (cf.~Urquhart  \cite[Lemma 6]{U78}).
If the lattice $\Lalg$ is distributive,
then $\mph{\Lalg}{\twoB}= \CL(\Lalg,\twoB)$,
the relation
$E$
coincides with the pointwise partial order
of maps and $\D{\Lalg}$ is the usual dual space
of
$\Lalg$
in the Priestley duality~\cite{P70}.

Plo\v s\v cica's representation
of $\Lalg \in \CL$
is then obtained in the following way.
One takes the family of continuous $E$-preserving partial maps
from $\D{\Lalg}$ into $\twoT_\T$, where $E$ on $\twoT$ is
taken to be the
usual order.
The maximally-defined members of this set are then shown to form a lattice
isomorphic to~$\Lalg$, the isomorphism being given by
the natural evaluation map $e_{\Lalg} \colon \Lalg \to
\mathrm{ED}(\Lalg)$.
We recall further details of the construction in Section~\ref{newCE}.

Our strategy for obtaining the canonical extension
$\Lalg^\delta$  of $\Lalg$
will be to replace $\D{\Lalg}$ above by $\Dflat{\Lalg}$. Here $\mathrm{D}^\flat$ is
the composition $^\flat \circ \mathrm{D}$, where $\flat$ is the
map
forgetting the topology.
Thus the category we shall use to build the canonical extension will be
the category  $\CG$
of graphs $\X=(X,E)$ and partial maps which preserve~$E$.
(We warn that what we are calling a \emph{graph} would usually be referred to as a \emph{digraph}.
We also warn
that  structures of the form
$\Dflat{\Lalg}$  will have special properties we shall need to
exploit in due course.)
Much of our work leading up to Theorem~\ref{canext} involves checking
that the ancillary results on
the graphs with topology
which underlie the Plo\v s\v cica
representation \cite[Section~1]{Pl95} have appropriate analogues in
the topology-free setting.

Initially we let $\X =(X,E)$ be any graph.
For two graphs $\X=(X,E_X)$ and $\Y=(Y,E_Y)$ we use the
notation $\CG(\X,\Y)$ to denote the collection of
total $E$-preserving maps from $\X$ to $\Y$. In Section \ref{cat}
we will make use of the set of partial morphisms,
$\CG^{\mathrm p}(\X,\Y)$.
We note that, by Zorn's Lemma, every partial $E$-preserving map into
$\twoT$
can be extended to a maximal partial $E$-preserving map;
by `maximal' we mean here that there is no partial $E$-preserving map properly extending it
(in general such an extension will not be unique).
We will denote by $\mpe{\X}{\twoT}$ the set of maximal partial
$E$-preserving maps from $\X$ to $\twoT$.
Let us use
the shorthand MPE to refer to an element of such a set.

In the distributive case we
can restrict to the situation in which
$E$
is a partial order and
MPE's
are simply total maps
which are order-preserving.  It is immediate that the order-preserving maps from a poset~$\X$
into~$\twoT$ always form a complete lattice under pointwise join and meet.
In the case of MPE's on a graph the situation is more complicated.

We now work towards showing, for a general graph $\X =(X,E)$, that
$\mpe{\X}{\twoT}$ does
indeed form a complete lattice.
As we shall see in Proposition~\ref{prop:compl}
and the discussion preceding it,  the similarity between Lemma~\ref{1.3-mod}
below and the corresponding result, \cite[Lemma 1.3]{Pl95}, concerning
graphs with topology,
is critical to the success of our approach to completions
built from dual spaces. The key point is that in \cite[Lemma 1.3]{Pl95}
the topological assertions and those involving~$E$ alone can be separated.

\begin{lem}\label{1.3-mod}
 {\upshape (cf. \cite[Lemma 1.3]{Pl95})}
Let\, $\X=(X,E)$ be a graph and
$\varphi \in \mpe{\X}{\twoT}$. Then
\begin{newlist}
\item[{\upshape (i)}]$\varphi^{-1}(0)=\{\, x\in X \mid \text{there is no}\ y\in\varphi^{-1}
(1)\ \text{with}\ (y,x)\in E \,\}$;
\item[{\upshape (ii)}] $\varphi^{-1}(1)=\{\, x\in X \mid
\text{there is no}\ y\in\varphi^{-1}
(0)\ \text{with}\ (x,y)\in E\,\}$.
\end{newlist}
\end{lem}

Consider
a family
$\{\, \varphi_i \in \mpe{\X}{\twoT} \mid i \in I \,\}$
for an index set $I$.
Then we can define
the
pointwise meet $p_{\varphi_i}^{\bigwedge}$ and pointwise
join $p_{\varphi_i}^{\bigvee}$ as follows:
\[
p_{\varphi_i}^{\bigwedge}(x)=
\begin{cases} 1&\text{if}\ x\in \bigcap_{i\in I} \varphi_i^{-1}(1),\\
             0&\text{if}\ x\in \bigcup_{i\in I} \varphi_i^{-1}(0)
\end{cases}  \text{ \quad and \quad } 
p_{\varphi_i}^{\bigvee}(x)=
\begin{cases} 1&\text{if}\ x\in \bigcup_{i\in I} \varphi_i^{-1}(1),\\
             0&\text{if}\ x\in \bigcap_{i\in I} \varphi_i^{-1}(0).
 \end{cases}
\]
(The motivation for these definitions comes from Plo\v s\v cica \cite[Section~3]{Pl95}.)
To see that these partial maps are $E$-preserving,
let $(p_{\varphi_i}^{\bigwedge}(x),p_{\varphi_i}^{\bigwedge}(y))=(1,0)$. Then
there exists
$j\in I$ with $\varphi_j(y)=0$,
and thus $(\varphi_j(x),\varphi_j(y))=(1,0)$
and so $(x,y)\notin E$.

We shall work with specific extensions of the maps $p_{\varphi_i}^{\bigwedge}$
and  $p_{\varphi_i}^{\bigvee}$ to elements of $\mpe{\X}{\twoT}$.
The  motivation for our choice
of maximal maps  comes
from Lemma~\ref{1.3-mod}.  Define
$e_{\varphi_i}^{\bigwedge} \colon \X\to \twoT$ and $e_{\varphi_i}^{\bigvee} \colon
\X\to \twoT$ as follows:
\begin{align*}
e_{\varphi_i}^{\bigwedge}(x)&=
\begin{cases} 1 &\text{if}\ x\in \bigcap_{i\in I} \varphi_i^{-1}(1),\\
             0 &\text{if}\ \text{there is no }
y\in\bigcap_{i\in I} \varphi_i^{-1}(1)\ \text{with }
(y,x)\in E;
\end{cases}  \\
e_{\varphi_i}^{\bigvee}(x)&=
\begin{cases} 1 &\text{if there is no}\ y\in\bigcap_{i\in I} \varphi_i^{-1}(0)\ \text{with}\ (x,y)\in E,\\
             0 &\text{if}\ x\in \bigcap_{i\in I} \varphi_i^{-1}(0).
\end{cases}
\end{align*}
We claim that $e_{\varphi_i}^{\bigwedge}$ extends $p_{\varphi_i}^{\bigwedge}$.
Let $x\in \bigcup_{i\in I} \varphi_i^{-1}(0)$, so that
$x\in \varphi_j^{-1}(0)$ for some $j\in I$. Then there is no
$y\in\bigcap_{i\in I} \varphi_i^{-1}(1)$ with $(y,x)\in E$, for otherwise
$(\varphi_j(y),\varphi_j(x))=(1,0)$, which contradicts the
preservation of~$E$ by $\varphi_j$.
Analogously,
$e_{\varphi_i}^{\bigvee}$ extends $p_{\varphi_i}^{\bigvee}$.
The next lemma shows that $e_{\varphi_i}^{\bigwedge}$ and
$e_{\varphi_i}^{\bigvee}$ are maximal partial $E$-preserving extensions of
$p_{\varphi_i}^{\bigwedge}$ and $p_{\varphi_i}^{\bigvee}$, respectively.

\begin{lem}\label{ext}
Let\, $\X=(X,E)$ be a graph  and let\, $\{\,\varphi_i \mid i\in I\,\} \subseteq \mpe{\X}{\twoT}$.
Then the maps\, $e_{\varphi_i}^{\bigwedge}$ and\, $e_{\varphi_i}^{\bigvee}$ are
elements of\, $\mpe{\X}{\twoT}$
extending\, $p_{\varphi_i}^{\bigwedge}$ and\, $p_{\varphi_i}^{\bigvee}$,
respectively.
\end{lem}

\begin{proof}
To see that $e_{\varphi_i}^{\bigwedge}$ preserves $E$,
consider $(e_{\varphi_i}^{\bigwedge}(x_1),e_{\varphi_i}^{\bigwedge}(x_0))=(1,0)$.
Then $x_1\in \bigcap_{i\in I} \varphi_i^{-1}(1)$ and
there is no $y\in \bigcap_{i\in I} \varphi_i^{-1}(1)$ with $(y,x_0)\in E$.
Thus $(x_1,x_0)\notin E$ as required.

We now
show the maximality of $e_{\varphi_i}^{\bigwedge}$.
Let $\psi \in \mpe{\X}{\twoT}$ be a map extending $e_{\varphi_i}^{\bigwedge}$.
Let $x\in X$ be such that $x\notin \dom(e_{\varphi_i}^{\bigwedge})$;
thus $x\notin \bigcap_{i\in I} \varphi_i^{-1}(1)$ and there is
$y_1\in \bigcap_{i\in I} \varphi_i^{-1}(1)$ with $(y_1,x)\in E$.
We want to show that $x\notin \dom(\psi)$.
Let $j\in I$ be such that $x\notin \varphi_j^{-1}(1)$.
Since $\varphi_j \in \mpe{\X}{\twoT}$, by Lemma~\ref{1.3-mod}
there exists $y_0\in \varphi_j^{-1}(0)$ with
$(x,y_0)\in E$. Hence $\varphi_j(y_0)=0$ and so,
by Lemma~\ref{1.3-mod} again, there is no $y\in  \varphi_j^{-1}(1)$
with $(y,y_0)\in E$. Hence
there is no $y\in  \bigcap_{i\in I} \varphi_i^{-1}(1)$ with
$(y,y_0)\in E$, which gives us $e_{\varphi_i}^{\bigwedge}(y_0)=0$.
As $\psi$ extends $e_{\varphi_i}^{\bigwedge}$, we obtain $\psi(y_0)=0$.
Moreover, as $e_{\varphi_i}^{\bigwedge}(y_1)=1$,  we have $\psi(y_1)=1$.

Now suppose for a contradiction that $x\in \dom(\psi)$.
Then $(y_1,x)\in E$ and $(x,y_0)\in E$ give us
$1=\psi(y_1)\le \psi(x) \le \psi(y_0)=0$, which is false. \qed
\end{proof}

\begin{thm}\label{compl}
Let\, $\X = (X,E)$ be a graph.
Then the set\, $\fC{\X}=\mpe{\X}{\twoT}$ ordered by the rule
\[
\varphi\le\psi \iff \varphi^{-1}(1)\subseteq\psi^{-1}(1)
\]
is a complete lattice.
\end{thm}

\begin{proof}
Obviously, the relation $\le$ is reflexive and transitive.
The antisymmetry of $\le$ follows from the fact that for any
$\varphi, \psi \in \mpe{\X}{\twoT}$, if
$\varphi^{-1}(1)=\psi^{-1}(1)$, then, by Lemma~\ref{1.3-mod},
$\varphi^{-1}(0)=\psi^{-1}(0)$, and hence $\varphi=\psi$.

To show that $\le$ is a (complete) lattice order, let
$\{\,\varphi_i \mid i\in I\,\} \subseteq \mpe{\X}{\twoT}$.
We claim that  $\bigwedge \varphi_i = e_{\varphi_i}^{\bigwedge}$
which is an element of $\fC{\X}$ by Lemma~\ref{ext}.

First note that
$(e_{\varphi_i}^{\bigwedge})^{-1}(1) = \bigcap_{i\in I} \varphi_i^{-1}(1) \subseteq \varphi_i^{-1}(1)$
for all $i\in I$.
Thus in $(\fC{\X},\le)$, we have
$e_{\varphi_i}^{\bigwedge} \le \varphi_i$ for all $i\in I$.
Now let
$\psi \in \mpe{\X}{\twoT}$ such that
$\psi \le \varphi_i$ for all $i\in I$. Then
$\psi^{-1}(1) \subseteq (\bigcap_{i\in I} \varphi_i)^{-1}(1) = (e_{\varphi_i}^{\bigwedge})^{-1}(1)$,
whence
$\psi \le e_{\varphi_i}^{\bigwedge}$ as required. \qed
\end{proof}

Before we move on to incorporate topology  we insert a
lemma
(the proof of which is straightforward)
concerning maps between graphs.
We shall need this in Section~\ref{cat}, when we consider morphisms acting by
composition.

\begin{lem} \label{E-pres} Let\, $\X=(X,E_X)$ and\, $\Y=(Y,E_Y)$ be graphs
and let\,
$\alpha \in \CG(\X,\Y)$.
Then
for\, $\varphi \in \CG^{\rm p}(\Y,\twoT)$,
the map\,
$\varphi \circ \alpha$, with
the domain being  determined by the set of elements for
which the composition is defined,
is a partial  map from\, $\X$ to $\twoT$ which preserves~$E$.
\end{lem}

A set equipped with a reflexive binary relation and a topology will be called a
\defn{graph with topology}.
We make the family
$\CG_\T$
of
graphs with topology
into a category in the following way.
A map $\varphi \colon (X_1,E_1,\tau_1)\to (X_2,E_2,\tau_2)$ between
graphs with topology
is called a
\defn{$\CG_\T$-morphism}
if it preserves
the binary relation and
is continuous as a map from  $(X_1,\tau_1)$ to  $(X_2,\tau_2)$.
A partial map
$\varphi\colon  (X_1,E_1,\tau_1)\to (X_2,E_2,\tau_2)$
is called a
\defn{partial $\CG_\T$-morphism}
if its domain is a $\tau_1$-closed subset of $X_1$ and the
restriction of $\varphi$ to its domain is a morphism. (We assume
that $\dom(\varphi)$ inherits the binary relation and the topology
from $X_1$.)
A partial
$\CG_\T$-morphism
is called maximal,
or an MPM for short,
if there is no partial
$\CG_\T$-morphism
properly extending it. For a
graph with topology
$\X_\T = (X,E,\T)$ we denote by
$\mpm{\X_\T}{\twoT_\T}$ the set of MPM's from $\X_\T$ to $\twoT_\T$.

Our candidate for the canonical extension of a bounded lattice~$\Lalg$
is the complete lattice $\mpe{\Dflat{\Lalg}}{\twoT}
=\fC{\Dflat{\Lalg}}$.
The embedding will be given, as in the Plo\v{s}\v{c}ica representation
(see Proposition~\ref{lem-eval} below) by
the map $e_{\Lalg}$, given by evaluation,
onto
the set of  maps which are maximal among continuous partial morphisms into $\twoT_\T$.
It is therefore necessary to reconcile the two versions of
maximality---one with topology and the other without.
This can be done for graphs in general, rather than just those
arising from dual spaces.

\begin{prop} \label{prop:reco}
Let\, $\X_\T = (X,E, \T)$ be a
graph with topology
and\, $\X = (X,E)$ be its untopologised counterpart.
Then\, $\mpm{\X_\T}{\twoT_\T} \subseteq \mpe{\X}{\twoT}$.
\end{prop}

\begin{proof}
We have already observed that
 Lemma~\ref{1.3-mod}
has a topological counterpart
shown in
\cite[Lemma~1.3]{Pl95}.
Specifically, for $\varphi \in \mpm{\X_\T}{\twoT_\T}$ we have
\begin{align*}
\varphi^{-1}(0) &=\{\, x\in X \mid \text{ there is no }
y\in\varphi^{-1} (1) \text{ with } (y,x)\in E \,\};\\
\varphi^{-1}(1) &=\{\, x\in X \mid
\text{ there is no } y\in\varphi^{-1}
(0) \text{ with } (x,y)\in E\,\}.
\end{align*}
Comparing
the above
with the properties of $\mpe{\X}{\twoT}$ in
Lemma~\ref{1.3-mod}
we assert that any $\varphi \in \mpm{\X_\T}{\twoT_\T}$ is maximal as a partial
 $E$-preserving map
when the topology is suppressed.
To verify this, suppose for a contradiction that  $\psi \in \mpe{\X}{\twoT}$
strictly extends $\varphi \in \mpm{\X_\T}{\twoT_\T}$.  Let $x \in \dom(\psi ) \setminus
\dom(\varphi)$ and assume without loss of generality that
$\psi (x) = 1$.
By maximality of $\psi$, Lemma~\ref{1.3-mod} implies that there
is no $y \in \dom(\psi ) $ such that $\psi (y) = 0$ and $(x,y) \in E$.
In particular, there is no $y \in \dom(\varphi) $ such that $\varphi(y) = 0$ and $(x,y) \in E$.
Now, by the characterisation of
$\varphi^{-1}(1)$ for $\varphi \in \mpm{\X_\T}{\twoT_\T}$, we deduce that $x \in  \varphi^{-1}(1)$,
which is the required contradiction. \qed
\end{proof}

The set $\mpm{\X_\T}{\twoT_\T}$ can
be considered as a subposet of the poset $\mpe{\X}{\twoT}$
with the partial order on $\mpm{\X_\T}{\twoT_\T}$ given, as in \cite{Pl95}, by
\[
\varphi\le\psi\quad\text{iff }\ \varphi^{-1}(1)\subseteq\psi^{-1}(1).
\]
If in particular $\X_\T = \D{\Lalg}$ then it was shown in \cite{Pl95} (see Proposition~\ref{lem-eval}(iii) below)
that the partial order $\le$ on $\mpm{\X_\T}{\twoT_\T}$ is a lattice order, and the lattice $(\mpm{\X_\T}{\twoT_\T},\le)$
is clearly a sublattice of the lattice
$\fC{\X} = (\mpe{\X}{\twoT},\le)$ (see Theorem~\ref{compl}).

\section{Canonical extensions constructed via Plo\v s\v cica's representation} \label{newCE}

Our next task is to establish that, when we take a bounded lattice~$\Lalg$
and the graph $\X = \Dflat{\Lalg}$, then $\fC{\X}$ acts as a completion of~$\Lalg$.
For this we call on the
Plo\v s\v cica representation,
and so now need to recall further facts from \cite{Pl95}.
The relationship between a lattice $\Lalg$ and the set
$\mpm{\D{\Lalg}}{\twoT_\T}$ is summed up in the following result.

\begin{prop}
{\upshape (\cite[Lemmas 1.2 and 1.5 and Theorem~1.7]{Pl95})}
\label{lem-eval}
Let\, $\Lalg\in \CL $
and let\, $\D{\Lalg}$ be the
graph with topology
defining the dual space of\,~$\Lalg$.
For\, $a\in L$, let the evaluation
map\, $e_a \colon \D{\Lalg} \to \twoT_\T$ be defined by
\[
e_a(f)= \begin{cases}
 f(a) &\text{ $a\in \dom(f)$,}\\
  - &\text{ undefined otherwise.}
    \end{cases}
       \]
Then the following hold.
\begin{newlist}
\item[{\upshape (i)}] The map\, $e_a \in \mpm{\D{\Lalg}}{\twoT_\T}$ for each\, $a \in L$.
\item[{\upshape (ii)}] Every\,
$\varphi \in \mpm{\D{\Lalg}}{\twoT_\T}$ is of the form\, $e_a$ for some\, $a\in L$.
\item[{\upshape (iii)}]
The map\, $e_{\Lalg}: \Lalg \to \mpm{\D{\Lalg}}{\twoT_\T}$ given by evaluation,\, $a \mapsto e_a$
{\upshape(}$a \in L${\upshape)}, is an
isomorphism of\, $\Lalg$ onto the lattice\, $\mpm{\D{\Lalg}}{\twoT_\T}$,
ordered by\, $\varphi\le\psi$ if and only if\, $\varphi^{-1}(1)\subseteq\psi^{-1}(1)$.
\end{newlist}
\end{prop}

By combining preceding results we
obtain the following
proposition.

\begin{prop} \label{prop:compl} Let\, $\Lalg \in\CL$ and let
\,$\D{\Lalg}$ be
the
graph with topology
dual to\,~$\Lalg$ and let \,$\X = \Dflat{\Lalg}$.
Then \,$(e, \fC{\X})$ is a
completion of\,~$\Lalg$, where \,$e \colon a \mapsto e_a$ {\upshape(}$a \in
L${\upshape)}.
\end{prop}

\begin{proof}
Theorem~\ref{compl} tells us that
$\fC{\X}$
is a complete lattice.
The result now follows
directly
from
Proposition~\ref{lem-eval}
combined with the fact that every $\varphi \in \mpm{\D{\Lalg}}{\twoT_\T}$,
and in particular
any
evaluation map $e_a$,
belongs to
$\mpe{\Dflat{\Lalg}}{\twoT}$, by
Proposition~\ref{prop:reco}. \qed
\end{proof}

 Proposition~\ref{prop:compl} identifies a completion for any  bounded lattice
constructed from the dual space of the lattice.  When
the lattice is distributive this certainly does give the canonical extension as
introduced by Gehrke and J\'onsson \cite{GJ94}.
We would now like to prove that this completion
supplies a canonical extension for an arbitrary bounded lattice.
To achieve this, we  need to
examine more closely
the structure of the dual space $\D{\Lalg}$ of a bounded lattice~$\Lalg$.

In \cite[Section~2]{Pl95}, Plo\v s\v cica demonstrates how his dual representation
for lattices relates
to
Urquhart's topological representation~\cite{U78}.
At  the level of the dual spaces, the passage back and forth between
Urquhart's
dual
representation and Plo\v s\v cica's
is set up
by  a bijection  between maximal
disjoint
filter-ideal pairs in $\Lalg$
(as employed by Urquhart)
and $\mph{\Lalg}{\twoB}$.
Instead of carrying a single binary relation~$E$,
Urquhart's dual spaces are equipped with a pair of quasi-orders,
$\le_1$ and $\le_2$.
Interpreted in terms of MPH's,
these  two relations are
defined on the set $\mph{\Lalg}{\twoB}$ as follows:
\[
f\le_1g \iff f^{-1}(1)\subseteq g^{-1}(1) \quad\text{and}\quad f\le_2g \iff
f^{-1}(0)\subseteq g^{-1}(0).
\]
These  quasi-orders $\le_1$ and $\le_2$ prove to be a valuable ancillary
tool for working with
graphs of the form
$\Dflat{\Lalg}=(\mph{\Lalg}{\twoB},E)$,
and we shall use them, in an
analogous
manner but on bigger domains,
in Section~\ref{cat}.

\begin{lem}
{\upshape (\cite[Theorem 2.1]{Pl95})}
\label{2.1}
Let\, $\Lalg\in \CL$  and let $f,g \in \mph{\Lalg}{\twoB}$. Then
\begin{newlist}
\item[{\upshape (i)}] $(f,g)\in E$ if and only if there exists $h \in \mph{\Lalg}{\twoB}$ with $f\le_1 h$ and $g\le_2 h$;
\item[{\upshape (ii)}] $f\le_2g$ if and only if there is no $h \in \mph{\Lalg}{\twoB}$ with $(h,g)\in E$ and
$(h,f)\notin E$;
\item[{\upshape (iii)}] $f\le_1g$ if and only if there is no $h \in \mph{\Lalg}{\twoB}$ with $(g,h)\in E$ and
$(f,h)\notin E$.
\end{newlist}
\end{lem}

It is a  consequence of (i) above
that $\le_1\,  \subseteq E$ and $\ge _2 \, \subseteq E$.
This will be used repeatedly below.

\begin{lem}\label{1.1-gen} Let\, $\Lalg\in \CL$.  Let\, $\X = \Dflat{\Lalg}$ and let\,
$\varphi \in \mpe{\X}{\twoT}$. Then
\begin{newlist}
\item[{\upshape (i)}] if $f \notin \varphi^{-1}(0)$ there exists
$g \in \varphi^{-1}(1)$ such that $f \le_2 g$;
\item[{\upshape (ii)}] if $f \notin \varphi^{-1}(1)$ there exists $ g \in \varphi^{-1}(0)$ such that $ f \le_1 g$.
\end{newlist}
\end{lem}

\begin{proof} (i) Given $f \notin \varphi^{-1}(0)$, we have by Lemma \ref{1.3-mod} that there exists $h \in \varphi^{-1}(1)$
such that $(h,f) \in E$. By Lemma \ref{2.1}(i) there exists
$g\in X$
such that $h \le_1 g$ and $f \le_2 g$, and we now claim that $g \in \varphi^{-1}(1)$. If we suppose that $g \notin \varphi^{-1}(1)$,
by Lemma \ref{1.3-mod} there must exist $u \in \varphi^{-1}(0)$ such that $(g,u) \in E$. Now $(g,u) \in E$ if and only if there exists
$v\in X$
such that $g \le_1 v$ and $u \le_2 v$, again by Lemma~\ref{2.1}(i). By the transitivity of $\le_1$ we have that
$h \le_1 v$.
Applying this with $u\le_2 v$ gives $(h,u) \in E$. Since $u \in \varphi^{-1}(0)$ and $h \in \varphi^{-1}(1)$, this contradicts that $\varphi$ is $E$-preserving. Thus $g \in \varphi^{-1}(1)$. The
 proof of~(ii) is similar. \qed
\end{proof}

\begin{lem}\label{cor1.1-gen} Let\, $\Lalg\in \CL$,
$\X = \Dflat{\Lalg}$ and \,$\varphi \in \mpe{\X}{\twoT}$.
For $f,g \in \mph{\Lalg}{\twoB}$
\begin{newlist}
\item[{\upshape (i)}] if $f \le_2 g$ and\, $\varphi(f) =0$, then\, $\varphi(g)=0$;
\item[{\upshape (ii)}] if $f \le_1 g$ and\, $\varphi(f) =1$, then\, $\varphi(g)=1$.
\end{newlist}
\end{lem}
\begin{proof} Let $f \le_2 g$ and $\varphi(f)=0$ and suppose that
$g \notin \varphi^{-1}(0)$. Then by Lemma~\ref{1.1-gen}(i) there exists
$h \in \varphi^{-1}(1)$ such that $g \le_2 h$. By the transitivity of $\le_2$,
we then see that $f \le_2 h$ and so $(h,f) \in E$.
Now $\varphi(h)=1 \nleqslant 0 = \varphi(f)$, contradicting
the fact that $\varphi$ is $E$-preserving. Thus $\varphi(g)=0$. \qed
\end{proof}

We can now show that $\mpe{\Dflat{\Lalg}}{\twoT}$, as an extension of $\mpm{\D{\Lalg}}{\twoT_\T}$, satisfies
the density condition. That is, every element of $\mpe{\Dflat{\Lalg}}{\twoT}$ can be written as
both a join of meets and a meet of joins of elements of $\mpm{\D{\Lalg}}{\twoT_\T}$.

\begin{prop} {\rm (}{\bf Density}{\rm )} \label{dens-direct}
Let\, $\Lalg\in \CL$  and let\, $\X = \Dflat{\Lalg}$. Then every
element\,
$\varphi$ of the completion \,$\fC{\X} = \mpe{\Dflat{\Lalg}}{\twoT}$
can be expressed as
$$
\varphi = \dJ \Big\{\, \DM F \bigm| \DM F \le \varphi \,\Big\}
= \DM \Big\{\, \dJ I \bigm| \varphi \le \dJ I \,\Big\}
$$
where $F$ ranges over the filters of\, $\mpm{\D{\Lalg}}{\twoT_\T}$ and~$I$ over the
ideals.
\end{prop}
\begin{proof} If
we let
$\psi = \bigvee \big\{\, \bigwedge F \mid  \bigwedge F \le \varphi \,\big\}$,
it is clear that $\psi \le \varphi$. In order to show that $\varphi \le \psi$, we must
check that
$\varphi^{-1}(1) \subseteq \psi^{-1}(1)$. Therefore, using the definitions of joins and meets
in $\mpe{\Dflat{\Lalg}}{\twoT}$ in Lemma~\ref{ext},
and given $x \in \varphi^{-1}(1)$, we must show that if
$y \in \bigcap\, \{\, (\bigwedge F)^{-1}(0) \mid \bigwedge F \le \varphi \,\} $,
then $(x,y) \notin E$.

In order to
do this, we claim that
$\bigcap\, \{\, (\bigwedge F)^{-1}(0) \mid \bigwedge F \le \varphi \,\} \subseteq \varphi^{-1}(0)$.
Suppose that $f \notin \varphi^{-1}(0)$. By
Lemma \ref{1.1-gen}(i) there exists
$g \in \varphi^{-1}(1)$ such that $f\le_2 g$. Now
we consider
$\mathcal{F}_g := \{\, a \in L \mid g(a)=1 \,\}$,
a filter of $\Lalg$.
From Proposition~\ref{lem-eval}(iii), we
can see that a set
$F \subseteq \mpm{\D{\Lalg}}{\twoT_\T}$
is a filter if and only if
$F = \{\, e_b \mid b \in \mathcal{F} \,\}$ where $\mathcal{F}$ is a filter of $L$.
This gives us a filter $G$
of
$\mpm{\D{\Lalg}}{\twoT_\T}$, with
$G = \{\, e_a \mid a \in \mathcal{F}_g \,\}$, and
we will show that $\bigwedge G \le \varphi$. If $h \in (\bigwedge G)^{-1}(1)$, then
$h \in \bigcap\, \{\, e_{a}^{-1}(1) \mid a \in \mathcal{F}_g \,\}$.
That is, for every $a \in L$ such that $g(a)=1$, $h(a)=1$.
So $g^{-1}(1) \subseteq h^{-1}(1)$, giving us that $g \le_1 h$.
Now we use
Lemma~\ref{cor1.1-gen}(ii) \textbf{}to get $\varphi(h)=1$ and hence $(\bigwedge G)^{-1}(1) \subseteq \varphi^{-1}(1)$.

We have $\bigwedge G \le \varphi$ and now show that $f \notin (\bigwedge G)^{-1}(0)$.
By definition
$$
\big(\, \DM G\, \big)^{-1}(0) = \Big\{\, b \in \mph{\Lalg}{\twoB} \bigm| \text{for all } c
\in \bigcap\, \{\, e_{a}^{-1}(1) \mid a \in \mathcal{F}_g \,\}, (c,b)\notin E \,\Big\}
$$
and clearly $g \in \bigcap\,\{\, e_{a}^{-1}(1) \mid a \in \mathcal{F}_g \,\}$.
Now,
our original assumption that
$f \le_2 g$ gives us that $(g,f) \in E$ and so
$f \notin (\bigwedge G)^{-1}(0)$.
Clearly now $f \notin \bigcap\, \{\, (\bigwedge F)^{-1}(0) \mid \bigwedge F \le \varphi \,\}$ and
thus we have
$\bigcap\, \{\, (\bigwedge F)^{-1}(0) \mid \bigwedge F \le \varphi \,\} \subseteq \varphi^{-1}(0)$.

Now recall that $x \in \varphi^{-1}(1)$ and let
$y \in \bigcap\, \{\, (\bigwedge F)^{-1}(0) \mid  \bigwedge F \le \varphi \,\}$.
Then $(x,y) \in E$ would imply $\varphi (x) = 1 \nleqslant 0 = \varphi(y)$, contradicting
$\varphi$ being $E$-preserving. Thus $(x,y) \notin E$ and we have $x \in \psi^{-1}(1)$ as required. \qed
\end{proof}

In Proposition~\ref{dens-direct} we proved very directly that
$\fC{\Dflat{\Lalg}}$ supplies a dense completion of~$\Lalg$.  We now want
to prove that this completion is compact.  This is true, but not entirely
straightforward to prove.  We shall go via a
somewhat
circuitous route,
obtaining on the way characterisations of filter and ideal elements
in $\fC{\Dflat{\Lalg}}$ which are of independent interest.

To
describe the filter
and ideal
elements of the extension $\fC{\Dflat{\Lalg}}$
we need to look more
closely at the topology on $\D{\Lalg}=(\mph{\Lalg}{\twoB},E,\T)$.
We recall that $\T$ has subbasic closed sets of the form
$V_a = \{\, f \in \mph{\Lalg}{\twoB} \mid f(a) = 0 \,\}$
and $W_a = \{\, f \in \mph{\Lalg}{\twoB} \mid f(a)=1 \,\}$
for $a \in L$.
We observe that from
Lemma~\ref{cor1.1-gen},
for any $\varphi \in \mpe{\Dflat{\Lalg}}{\twoT}$,
we have that $\varphi^{-1}(1)$ is a $\le_1$-increasing set, and
$\varphi^{-1}(0)$ is a $\le_2$-increasing set. We also note that the intersection of a collection
of closed sets of the form $W_b$ for $b \in L$ will be
$\le_1$-increasing,
while the intersection of sets
of the form $V_a$ with $a \in L$ will be $\le_2$-increasing.

The following result parallels ~\cite[Lemma 1.4]{Pl95}
though our approach in proving it is different.

\begin{lem}\label{1.4-mod}
Let\, $\Lalg\in \CL$ and let
$a\in L$. Let\,
$\X =
{\Dflat{\Lalg}}$
and let\, $\varphi\in \mpe{\Dflat{\Lalg}}{\twoT}$. Then
\begin{newlist}
\item[{\upshape (i)}] $\varphi^{-1}(0)\cap W_a =\emptyset$ implies\, $\varphi^{-1}(0)\subseteq V_a$;
\item[{\upshape (ii)}] $\varphi^{-1}(1)\cap V_a =\emptyset$ implies\, $\varphi^{-1}(1)\subseteq W_a$.
\end{newlist}
\end{lem}

\begin{proof}
We consider (i).
Let $\varphi^{-1}(0)\cap W_a =\emptyset$. Suppose that
$f\in \varphi^{-1}(0)$ and $f\notin V_a$. As $V_a=e_a^{-1}(0)$, by Lemma~\ref{1.1-gen}(i) there exists $ g \in \varphi^{-1}(1)$ such that
$f \le_2 g$. Since $\varphi(f) =0$, by Lemma~\ref{cor1.1-gen}(i) we obtain $\varphi(g)=0$, a contradiction. \qed
\end{proof}

\begin{lem} \label{closedpreimages}
Let\, $\Lalg\in \CL$, let\, $\X_\T = \D{\Lalg}$ and
let\, $\varphi $ be an element of\,
$\mpe{\Dflat{\Lalg}}{\twoT}$. Then
\begin{newlist}
\item[{\upshape (i)}] $\varphi^{-1}(1)$ is  a\, $\T$-closed subset of\, $\X_\T$ if and only if
\,$\varphi^{-1}(1)= \bigcap\, \{\, W_b \mid b \in K \,\}$, for some
$K \subseteq L$;
\item[{\upshape (ii)}] $\varphi^{-1}(0)$ is a\, $\T$-closed subset of\, $\X_\T$ if and only if
\,$\varphi^{-1}(0)= \bigcap\, \{\, V_a \mid a \in M \,\}$,
for some $M \subseteq L$.
\end{newlist}
\end{lem}
\begin{proof}
We consider (i) and so let $X$ denote the underlying set of $\X_\T$,
that is,
$X = \mph{\Lalg}{\twoB}$.
The sufficiency
is clear as each $W_b$ is $\T$-closed. For the
necessity, let $\varphi^{-1}(1)$ be a $\T$-closed subset of $\X_\T$. Let
$$
\cat S := \{\,W_b \mid b\in L, \varphi^{-1}(1) \subseteq W_b\,\}
$$
be the
family of all sets $W_b$ containing $\varphi^{-1}(1)$ and
$K := \{\,b\in L \mid W_b\in \cat S
\,\}$.
As
 $\varphi^{-1}(1) \subseteq W_b$ for each $b\in K$, it is obvious that
$\varphi^{-1}(1)
\subseteq \bigcap\, \{\, W_b \mid b\in K\,\}$.
The reverse
inclusion will hold too provided we can show
that $f\notin \varphi^{-1}(1)$
implies
the
existence of $W_b\in \cat S$
such that $f\notin W_b$.

Let $f\notin \varphi^{-1}(1)$. By Lemma~\ref{1.3-mod}(ii), there
exists $g\in \varphi^{-1}(0)$ such that $(f,g)\in E$.
Since
$g\in
\varphi^{-1}(0)$, for every $h\in \varphi^{-1}(1)$ we necessarily have
$(h,g)\notin E$. By definition of the relation~$E$ this means that
there exists $b_h\in L$ such that $h(b_h)=1$ and $g(b_h)=0$. We note
that each set $X\setminus V_{b_h}$ is $\T$-open and $h\in
X\setminus V_{b_h}$. Hence
$\varphi^{-1}(1)\subseteq \bigcup\,\{\,X\setminus V_{b_h} \mid h\in \varphi^{-1}(1)\,\}$.
Since the space $\X_\T$
is compact and $\varphi^{-1}(1)$ is a $\T$-closed subset
of $\X_\T$ by hypothesis, there are $b_1,\dots,b_n\in L$ such that
$g(b_1)=\dots =g(b_n)=0$ and $\varphi^{-1}(1)\subseteq (X\setminus
V_{b_1}) \cup \dots \cup (X\setminus V_{b_n})$. Let $b:= b_1
\lor\dots\lor b_n\in L$. We have
$g(b)=0$
and $\varphi^{-1}(1)\subseteq
X\setminus V_b$. The latter by Lemma~\ref{1.4-mod} yields
$\varphi^{-1}(1)\subseteq W_b$,
and thus $W_b \in \cat S$.
Since $(f,g)\in E$ and
$g(b)=0$, we cannot have $f(b)=1$.
Hence $f\notin W_b$,
as required. \qed
\end{proof}

\begin{prop}\label{cl-op}  {\rm (}{\bf Filter
and ideal
elements}{\rm)}
Let\, $\Lalg\in \CL$,
$\X_\T = \D{\Lalg}$ and let\,
$\varphi$ be an element of \,$\mpe{\Dflat{\Lalg}}{\twoT} = \fC{\X}$.
Then the three conditions in {\upshape (i)} are
 equivalent and the three conditions in {\upshape (ii)} are equivalent:
\begin{newlist}
\item[{\upshape (i)}]
 \begin{itemize}
 \item[{\upshape (1)}] $\varphi$ is a filter
element of\, $\fC{\X}$;
 \item[{\upshape (2)}] $\varphi^{-1}(1)= \bigcap\, \{\, W_b \mid b \in K \,\}$, for some $K \subseteq L$;
 \item[{\upshape (3)}] $\varphi^{-1}(1)$ is a\, $\T$-closed subset of\, $\X_\T$.
 \end{itemize}
\item[{\upshape(ii)}]
 \begin{itemize}
 \item[{\upshape (1)}] $\varphi$ is an ideal
element of\, $\fC{\X}$;
 \item[{\upshape (2)}] $\varphi^{-1}(0)= \bigcap\, \{\, V_a \mid a \in M \,\}$, for some $M \subseteq L$;
 \item[{\upshape (3)}] $\varphi^{-1}(0)$ is a\, $\T$-closed subset of\, $\X_\T$.
 \end{itemize}
\end{newlist}
\end{prop}

\begin{proof}
The equivalences of (2) and (3) in~(i) and~(ii) come from Lemma~\ref{closedpreimages}.
Now let~(1) in~(i) hold, that is,
let
$\varphi$ be a filter element of $\fC{\X}$. Then
we have that
$\varphi = \bigwedge \psi_i$ where
$\{\,\psi_i \mid i\in I\,\} \subseteq \mpm{\D{\Lalg}}{\twoT_\T}$ for some index set $I$.
From the representation theorem
(see Proposition~\ref{lem-eval}),
 we know that each $\psi_i$ is in fact
$e_{b_i}$ for some $b_i \in L$ and from the definition of
a meet
of
elements of $\mpe{\Dflat{\Lalg}}{\twoT}$,
we have
$$ \varphi^{-1}(1) = \bigcap_{i \in I} \psi_i^{-1}(1) = \bigcap_{i \in I} e_{b_i}^{-1}(1)  = \bigcap_{i \in I} W_{b_i}.$$
Hence (1) implies (2). Now let (2) hold. Then
$\varphi^{-1}(1) = \bigcap_{i \in I} e_{b_i}^{-1}(1)$,
and thus $\varphi^{-1}(1)$
can be expressed as the
meet of images of lattice elements in $\fC{\X}$,
yielding (1).
The equivalence of~(1) and~(2) in~(ii) can be shown analogously. \qed
\end{proof}

We are now able to prove that our completion is indeed a compact completion.

\begin{prop} {\rm (}{\bf Compactness}{\rm )} \label{compact}
Let\, $\Lalg\in \CL$  and let\, $\X_\T = (\mph{\Lalg}{\twoB},E,\T)$. The lattice\,
$\fC{\X}=\mpe{\Dflat{\Lalg}}{\twoT}$
is a compact completion of\, $\Lalg$.
\end{prop}

\begin{proof} 
Let $A, B \subseteq L$ such that $\bigwedge \{\, e_a \mid a \in A\,\} \le \bigvee \{\, e_b \mid b \in B \,\}$ 
in the completion $\fC{\X}$.
Let $\varphi = \bigwedge \{\,e_a \mid a \in A \,\}$ and $\psi = \bigvee \{\, e_b \mid  b \in B \,\}$. We have
$$
\varphi^{-1}(1)= \bigcap_{a \in A} e_a^{-1}(1) = \bigcap_{a \in A} W_a\quad \text{and}\quad
\psi^{-1}(0)=
\bigcap_{b \in B} e_b^{-1}(0)
= \bigcap_{b \in B} V_b.
$$
Now $\varphi \le \psi$ 
in $\fC{\X}$ 
means $\varphi^{-1}(1)\subseteq \psi^{-1}(1) \subseteq X\setminus \psi^{-1}(0)$,
and hence
$$
\bigcap_{a \in A} W_a \subseteq X\setminus \bigcap_{b \in B} V_b
=\bigcup_{b \in B}  
\,(X\setminus V_b).
$$
Since the space 
$\X_\T$ is compact,
and since
each $W_a$ is $\T$-closed and each $X\setminus V_b$ is $\T$-open,
the topological compactness yields the existence of finite subsets
$A':=\{a_1,\dots,a_r\}\subseteq A$ and
$B':=\{b_1,\dots,b_s\}\subseteq B$ such that
$$W_{a_1}\cap \dots \cap W_{a_r} \subseteq (X\setminus V_{b_1}) \cup \dots \cup (X\setminus V_{b_s}).$$
We then let
$a':= a_{1} \land \dots \land a_{r}$ and $b':= b_{1}\lor \dots \lor b_{s}$ and observe that
$$
e_{a'}^{-1}(1) = W_{a'} = \bigcap_{a \in A'}
W_a \subseteq
\bigcup_{b \in B'}
(X\setminus V_b) = X\setminus V_{b'} = X \setminus e_{b'}^{-1}(0).$$

Now by Lemma~\ref{1.4-mod}(ii) we get $e_{a'}^{-1}(1) \subseteq e_{b'}^{-1}(1)$ and so
$\bigwedge \{\, e_a \mid a \in A'\,\} \le \bigvee \{\, e_b \mid b \in B' \,\}$ as required.
\qed 
\end{proof}

The principal
result of this section
is now an immediate  consequence of
Propositions~\ref{dens-direct} and \ref{compact}.

\begin{thm} {\rm (}{\bf Canonical extension}{\rm )} \label{canext}
Let\, $\Lalg\in \CL$ and let\, $\X = (\mph{\Lalg}{\twoB},E)$.
The lattice\, $\fC{\X}=\mpe{\X}{\twoT}$
ordered by
\[
\varphi\le\psi \iff \varphi^{-1}(1)\subseteq\psi^{-1}(1)
\]
is the canonical extension of\,
$\Lalg$.
\end{thm}

We
summarise
in Fig.~\ref{fig:Psquare} what we have achieved.
Here
\begin{alignat*}{2}
   & \mathrm D \colon \Lalg  &&\longmapsto  \X_\T := (\mph{\Lalg}{\twoB}, E, \T) , \\
   &\,^\flat\, \colon  \X_\T   && \longmapsto  \X := (\mph{\Lalg}{\twoB}, E), \\
   & \mathrm G \colon \X &&\longmapsto  \fC{\X}:=\mpe{\X}{\twoT},
   \end{alignat*}
 and   $\Lalg \hookrightarrow \Lalg^\delta \text{ via  } a \mapsto e_a$.
The fact that the diagram commutes is the content of
Theorem~\ref{canext}.

\begin{figure}  [ht]
\begin{diagram}[labelstyle=\textstyle]
\CL&
\rTo^{{\mathrm D}}
& \CG_\T\\
\dTo^{^\delta} & & \dTo_{\flat}\\
\CL^+&  \lTo_{{\mathrm G}} & \CG
\end{diagram}
\caption{Factorisation of $^\delta $ on objects
in $\CL$
 \label{fig:Psquare}}
\end{figure}

\section{A categorical framework for the canonical extension} \label{cat}

The canonical extension construction on $\CL$ is functorial; see
\cite{GH01,GP08}.  And, as we observed in
Section~\ref{intro},
there is, in the distributive case, a factorisation of
the canonical extension functor in terms of functors
obtained from the
hom-functors setting up Priestley duality.
In this section we seek an analogous result for~$\CL$,
insofar as this is possible.
Specifically we seek to set up a commutative diagram as shown in
Fig.~\ref{fig:Lsquare}, paralleling that shown for~$\CCD$ in
Fig.~\ref{fig:Dsquare}.

\begin{figure}  [ht]
\begin{diagram}[labelstyle=\textstyle]
\CL& \pile{
\rTo^{{\overline{\mathrm{D}}}}\\ \lTo_{{\overline{\mathrm{E}}}}} & \CY_\Tp \\
\dTo^{^\delta} & & \dTo_{\flat}\\
\CL^+&  \lTo_{{\overline{\mathrm{G}}}} & \CY
\end{diagram}
\caption{Seeking to factorise the functor~$^\delta$ on~$\CL$\label{fig:Lsquare}}
\end{figure}

In Fig.~\ref{fig:Lsquare}
we require
\begin{itemize}
\item $\CL^+$ to be the category of complete lattices with complete lattice homomorphisms,
\item $\CY_\Tp$ to be a category of graphs with topology and $\CY$ the corresponding category of graphs, obtained by forgetting the topology;
\item $\overline{{\mathrm{D}}}$ and $\overline{\mathrm{E}} $ to set up  a dual adjunction such that, for all $\Lalg \in \CL$, we have
$\Lalg \cong \overline{\mathrm{E}}\hspace*{.04cm}\overline{\mathrm{D}}(\Lalg)$
with the isomorphism given by a natural evaluation map,
\item $\overline{{\mathrm{D}}}$, $\overline{\mathrm{E}} $  and $\overline{\mathrm{G}}$ to be contravariant
functors,  sending   an object~$C$ into
some  appropriately specified set of partial maps into, respectively, $\twoB$, $\twoT_\T$ and~$\twoT$, and
\item  each of  $\overline{{\mathrm{D}}}$, $\overline{\mathrm{E}} $
 and
$\overline{\mathrm{G}}$  to act on morphisms by composition.
\end{itemize}

Thinking in terms of maps between objects, rather than functors,
we have seen (see Fig.~\ref{fig:Psquare}) that we can obtain
a diagram of this kind, based on the Plo\v s\v cica representation.
But regrettably
there are inherent obstacles to adding morphisms within
Plo\v s\v cica's framework.
 For $u \in
\CL(\Lalg, \K)$,
the morphism $\Dflat{u} \colon \Dflat{\K} \to \Dflat{\Lalg}$
would need to be given by $\Dflat{u}(f) := f \circ u$,
for all $f\in \Dflat{\K}$, where the domain is determined by the
set of elements for which the composition is
defined.   Certainly, for $f$ a maximal partial homomorphism,  $f \circ u$
is a partial homomorphism.
But, as the example in  Fig.~\ref{fig:nonmax} demonstrates,
$f \circ u$ need not be maximal.

We have already noted that Plo\v s\v cica's representation is a recasting of
Urquhart's topological  representation for bounded lattices.
Urquhart was able to set up a dual representation for surjective
$\CL$-morphisms but not for $\CL$-morphisms in general.
Therefore  our example  should come as no surprise.  Moreover,
it is very easy to see that $f \circ u$ is maximal if $u$ is
a surjective morphism and~$f$ an MPH.

A solution to the problem of making duality for lattices functorial
was proffered by Hartung \cite{GHext}, using
topological contexts, but this approach does not fit our  needs.
Our remedy is quite different.
We shall call on the topological representation of bounded lattices due to
Allwein and Hartonas~\cite{A92,AH93}.
This  is in the same style as that of Urquhart,
in that it makes use of
disjoint filter-ideal pairs of~$\Lalg$ to construct a dual space
for~$\Lalg$.
However such pairs are not required to be maximal, as they are in
Urquhart's representation.  Allwein and Hartonas thereby overcome the
problem encountered by Urquhart and are
able to set up a dual equivalence between~$\CL$ and a specified category of
topological structures.
This is achieved at the cost of working with a greatly
enlarged dual space.  A key feature of Priestley
duality  is its  `logarithmic' property (cf. ~\cite[Section~6.3]{CD98}).
For a finite lattice,
moving
to a dual space which is bigger than the lattice itself
constitutes  a major sacrifice,  and much of the appeal and power of
duality for distributive lattices is lost.
This doubtless   explains why the work of Allwein and
Hartonas has attracted little attention.
However the Allwein--Hartonas
duality, recast, insofar as this is necessary,
in terms of suitable partially defined morphisms,
provides exactly the formalism  we need  to achieve our aim of
linking the canonical extension construction to duality in a functorial way,
as in Fig~\ref{fig:Lsquare}.
We shall now outline what we have to do to carry out our programme.

\begin{figure}
\begin{center}
\begin{tikzpicture}[scale=0.7]
\po{1,1}
\po{1,0}
\po{1,-1}

\li{(1,-1)--(1,0)--(1,1)}

\node at (0.6,1) {$1$};
\node at (0.6,0) {$a$};
\node at (0.6,-1) {$0$};

\node at (1,-2) {$\Lalg$};

\path[->] (1.3,1) edge node[above] {$u$} (4.7,1);
\path[->] (1.3,0) edge node[above] {} (3.7,0);
\path[->] (1.3,-1) edge node[above] {} (4.7,-1);

\po{4,0}
\po{5,0}
\po{6,0}
\po{5,1}
\po{5,-1}

\li{(5,-1)--(4,0)--(5,1)}
\li{(5,-1)--(5,0)--(5,1)}
\li{(5,-1)--(6,0)--(5,1)}

\node at (4,-0.3) {$b$};
\node at (4.8,-0.3) {$c$};
\node at (6,-0.3) {$d$};
\node at (5,-2) {$\K$};

\path[->] (5.2,1) edge node[above] {$f$} (8.8,0.6);
\path[->] (5.1,0.1) edge[bend left=10] node[above] {} (8.8,0.4);
\path[->] (6.1,0.1) edge[bend left=10] node[above] {} (8.8,-0.4);
\path[->] (5.2,-1) edge node[above] {} (8.8,-0.6);

\po{9,-0.5}
\po{9,0.5}
\li{(9,-0.5)--(9,0.5)}

\node at (9.4,0.5) {$1$};
\node at (9.4,-0.5) {$0$};
\node at (9,-2) {$\twoB$};

\end{tikzpicture}
\caption{Failure of morphisms to act on MPH's by composition \label{fig:nonmax}}
\end{center}
\end{figure}

First of all, we require a Plo\v s\v cica-style presentation of the
Allwein--Hartonas dual spaces.
The first task is to recognise which
partial homomorphisms from~$\Lalg$ into~$\twoB$ correspond to
disjoint
filter-ideal
pairs
and appropriately to equip  the resulting set of partial homomorphisms with the
structure of a
graph with topology.
We may consider the full subcategory $\CY_\T$ of $\CG_\T$
whose objects are
the enlarged
dual spaces (in our graph formulation, given below) of members of~$\CL$.
 We  then let $\CY$ be the image of
$\CY_\T$ under~$^\flat$.
Because of their provenance,   we shall refer to
the objects of~$\CY$ as \defn{$\CL$-graphs}.

The results of Allwein and Hartonas \cite{A92,AH93}
imply that there are
contravariant functors
$\overline{\mathrm{D}}\colon
\CL \to \CY_\Tp$
and $\overline{\mathrm{E}} \colon \CY_\Tp \to \CL$ which are such that
$\Lalg \cong \overline{\mathrm{E}}\hspace*{.04cm}\overline{\mathrm{D}}(\Lalg)$
for every $\Lalg \in \CL$; we show that  the isomorphism is set  up by
evaluation
maps in the same manner as in
Proposition~\ref{lem-eval}
(see Propositions~\ref{eabar} and~\ref{mpm=ea}).
There is now   an obvious candidate for the canonical extension
of~$\Lalg$: we forget the topology on $\overline{\mathrm{D}}(\Lalg)$ and
take the set of all maximal partial morphisms into~$\twoT$.
 Rather than
verifying directly that we thereby obtain a dense and compact completion
we instead set up an order-isomorphism between our new candidate for the canonical extension
and the one based on the Plo\v s\v cica representation.  Finally we check out
that we really do have a well-defined functor
$\overline{\mathrm{G}}$ making the diagram in
Fig.~\ref{fig:Lsquare} commute, and at the same time use our diagram to confirm, in a transparent way,
that $^\delta $ does lift each $\CL$-morphism to a complete lattice homomorphism.
We now implement the strategy we have set out.

Let $\Lalg \in \CL$.  Then the  disjoint filter-ideal pairs of $\Lalg$
are obviously in bijective correspondence with those  partial homomorphisms $f$ into
$\twoB$,  for which
$(f^{-1}(1),f^{-1}(0))$ is a
(disjoint)
filter-ideal pair of $\Lalg$;
let us call them \defn{special partial homomorphisms},
SPH's for short,
and denote
the set of all such  maps by $\sph{\Lalg}{\twoB}$.
We note that every maximal partial homomorphism into $\twoB$ is special
but that it is easy to find examples of partial homomorphisms into $\twoB$ which fail to be special.
The presentation  by Allwein and Hartonas now permits a translation into
an equivalent formulation in terms of
special
partial homomorphisms (for the first dual) and continuous partial
morphisms (for the second dual).
The topology in \cite{AH93} is defined in the same way as for the
representation theorem
of Plo\v{s}\v{c}ica. That is, let sets of the form
\[
V_a=\{\,f \in \sph{\Lalg}{\twoB} \mid f(a)=0 \,\} \quad \text{and} \quad
W_a=\{\,f \in \sph{\Lalg}{\twoB} \mid f(a)=1 \,\},
\]
form a subbasis for the closed sets of the topology $\T$.
 We now define
the relation~$E$ on $\sph{\Lalg}{\twoB}$ by the rule
\[
(f,g)\in E  \text{ iff } f(x)\le g(x) \text{ for every }
x\in\dom(f)\cap\dom(g),
\]
that is, we take  the obvious extension of $E$, as defined  earlier
 on the subset
 $\mph{\Lalg}{\twoB}$. We
denote the dual space of $\Lalg$ by $\overline{{\mathrm{D}}}(\Lalg) =
(\sph{\Lalg}{\twoB}, E, \T)$.
The dual space specified in \cite{AH93}, called an \defn{enhanced L-space},
takes the form of  a topological space equipped with two quasi-orders ($\le_1$ and $\le_2$)
and a specified subbasis of closed sets.

We shall need the following two lemmas concerning
SPH's.
The first of these relates the relation~$E$ to quasi-orders
$\le_1$ and $\le_2$
(defined on $\sph{\Lalg}{\twoB}$ exactly as on $\mph{\Lalg}{\twoB}$).
For the case of $\mph{\Lalg}{\twoB}$
it appears as
Lemma~\ref{2.1}, which we carried over from \cite{Pl95};
here
we show that
the  result extends to SPH's in general.

\begin{lem}\label{2.1new}
Let\, $\Lalg\in \CL$  and let $f,g \in \sph{\Lalg}{\twoB}$. Then
\begin{newlist}
\item[{\upshape (i)}] $(f,g)\in E$ if and only if there exists $h \in \sph{\Lalg}{\twoB}$ with $f\le_1 h$ and $g\le_2 h$;
\item[{\upshape (ii)}] $f\le_2g$ if and only if there is no $h \in \sph{\Lalg}{\twoB}$ with $(h,g)\in E$ and
$(h,f)\notin E$;
\item[{\upshape (iii)}] $f\le_1g$ if and only if there is no $h \in \sph{\Lalg}{\twoB}$ with $(g,h)\in E$ and
$(f,h)\notin E$.
\end{newlist}
\end{lem}
\begin{proof}(i) If $(f,g) \in E$ then $f^{-1}(1) \cap g^{-1}(0)=\emptyset$.
Consider $h \in \sph{\Lalg}{\twoB}$ defined by
letting
$h^{-1}(1)=f^{-1}(1)$ and $h^{-1}(0)=g^{-1}(0)$.
Then $f \le_1 h$ and $g \le_2 h$. Conversely, if there exists $h \in \sph{\Lalg}{\twoB}$
with $f^{-1}(1) \subseteq h^{-1}(1)$ and $g^{-1}(0) \subseteq h^{-1}(0)$
then clearly $f^{-1}(1) \cap g^{-1}(0) = \emptyset$ and so $(f,g) \in E$.

(ii) Suppose that $f \le_2 g$ and let $h \in \sph{\Lalg}{\twoB}$.
Now suppose that $(h,f) \notin E$. This implies that there exists
$a \in L$ such that $h(a)=1$ and $f(a)=0$. Clearly $g(a)=0$
and so $(h,g) \notin E$.

Next, assume that
$f \nleqslant_2 g$. This implies that there exists $a \in L$ such that
$f(a)=0$ and $g(a) \neq 0$. Consider $h \in \sph{\Lalg}{\twoB}$ defined
by $h^{-1}(1)={\uparrow}a$ and $h^{-1}(0)=g^{-1}(0)$. Clearly $(h,g) \in E$
and $(h,f) \notin E$, completing the proof.
Part (iii) is proved in the same way. \qed
\end{proof}

We now derive  an
extension of Lemmas~\ref{cor1.1-gen}
and~\ref{1.4-mod} from $\Dflat{\Lalg}$ to
$\smash{\overline{\mathrm{D}}}^\flat(\Lalg)$:

\begin{lem}\label{3.5+3.7new} Let\, $\Lalg\in \CL$.
Let\, $\Y = (\sph{\Lalg}{\twoB},E)$ and let\,
$\varphi \in \mpe{\Y}{\twoT}$ and $f,g \in Y$. Then
\begin{newlist}
\item[{\upshape (i)}] if $f \le_2 g$ and\, $\varphi(f) =0$, then\, $\varphi(g)=0$;
\item[{\upshape (ii)}] if $f \le_1 g$ and\, $\varphi(f) =1$, then\, $\varphi(g)=1$.
\end{newlist}
Further, for $a \in L$, we have
\begin{newlist}
\item[{\upshape (iii)}] $\varphi^{-1}(0)\cap W_a =\emptyset$ implies\, $\varphi^{-1}(0)\subseteq V_a$;
\item[{\upshape (iv)}] $\varphi^{-1}(1)\cap V_a =\emptyset$ implies\, $\varphi^{-1}(1)\subseteq W_a$.
\end{newlist}
\end{lem}

\begin{proof}
Let $f \le_2 g$ and $\varphi(f) =0$. Then by Lemma~\ref{1.3-mod}, for all
$h \in \varphi^{-1}(1)$ we have $(h,f) \notin E_Y$, that is,
$h^{-1}(1) \cap f^{-1}(0) \neq \emptyset$. From $f \le_2 g$ we have
$f^{-1}(0) \subseteq g^{-1}(0)$, and then clearly for all
$h \in \varphi^{-1}(1)$ we have $h^{-1}(1) \cap g^{-1}(0) \neq \emptyset$,
whence $(h,g) \notin E_Y$. Thus $g \in \varphi^{-1}(0)$, again by Lemma~\ref{1.3-mod}.

Now consider
part (iii). Let $\varphi^{-1}(0)\cap W_a =\emptyset$ and
suppose that $f \notin V_a$ but $f \in \varphi^{-1}(0)$.
Since $f(a)\neq 0$, consider the
SPH
$g_f \colon \Lalg \to \twoB$ defined by $g_f^{-1}(1) = {\uparrow}a$ and
$g_f^{-1}(0)=f^{-1}(0)$. Now $g_f \in W_a$, but since $f \le_2 g_f$ we have
by~(i)
that $g_f \in \varphi^{-1}(0)$. \qed
\end{proof}

We are now ready to confirm that $\Lalg$ can be recaptured from the
graph with topology
$\overline{{\mathrm{D}}}(\Lalg)$
in the same manner
as, by Plo\v s\v cica's representation, it is recaptured from
$\D{\Lalg}$.
We first reformulate the representation due to
Allwein and Hartonas~\cite{A92,AH93}  in terms
of
graphs with topology
in the same way that Plo\v{s}\v{c}ica reformulated the~$L$-spaces of Urquhart.
We want to represent the lattice $\Lalg$ by the maximal partial continuous $E$-preserving
maps (MPM's) from $\overline{{\mathrm{D}}}(\Lalg)$ into $\twoT_\T$.

Let $\Lalg\in \CL $,
let
$\Y_\T = \overline{{\mathrm{D}}}(\Lalg)$ and $\Y=(\sph{\Lalg}{\twoB},E)$.
Let $a \in L$ and define
$\bar{e}_a \colon \Y_\T \to \twoT_\T$ for $f \in \sph{\Lalg}{\twoB}$ by
\[
\bar{e}_a(f) = \begin{cases}
f(a) & \text{ if } a \in \dom (f),\\
- & \text{ otherwise.}
\end{cases}
\]

\begin{prop} \label{eabar}   Let\, $\Lalg\in \CL$  and
define\, $\Y_\T = \overline{{\mathrm{D}}}(\Lalg)$, $\Y=(\sph{\Lalg}{\twoB},E)$, and\,
$\bar{e}_a$ as above.
Then\, $\bar{e}_a \in \mpm{\Y_\T}{\twoT_\T}$.
Furthermore, if\, $\bar{e}_a$  is regarded  as a map from
\,$\Y$ to $\twoT$, then
$\,\bar{e}_a~\in~\mpe{\Y}{\twoT}$.
\end{prop}

\begin{proof}
We note that $\bar{e}_a^{-1}(1)=W_a$ and $\bar{e}_a^{-1}(0)=V_a$
are closed sets and since the set
$\dom (\bar{e}_a) =W_a\cup V_a$ is closed, too,
we have that
$\bar{e}_a \colon \dom(\bar{e}_a) \to \twoT_\T$
is continuous.
Suppose $\bar{e}_a(f)=1$ and $\bar{e}_a(g)=0$. Then $(f,g) \notin E$ and so
$\bar{e}_a$ is $E$-preserving.
Now assume that $\phi \in \mpm{\Y_\T}{\twoT_\T}$ and the domain of $\phi \colon \Y_\T \to \twoT_\T$
properly extends the domain of $\bar{e}_a$. The first case to consider is
$f \in \dom (\phi) \setminus \dom (\bar{e}_a)$
such that $\phi(f)=1$. This implies that for all
$g \in \bar{e}_a^{-1}(0)$,
$(f,g) \notin E$.
In particular it means
that for $h \in \sph{\Lalg}{\twoT}$ with
$h^{-1}(1)=\{1\}$ and  $h^{-1}(0) = {\downarrow}a$ we have
$(f,h) \notin E$.
Since
$f^{-1}(1)$ is a filter, this  can only happen if
$f(a)=1$ and so
$f\in \dom (\bar{e}_a)$,
a contradiction.
Similarly one can show that if
$\phi(f)=0$ then
$f(a)=0$, a contradiction.
Thus the domain of $\bar{e}_a$ is maximal.

The fact that $\bar{e}_a \in \mpe{\Y}{\twoT}$ follows from
Proposition~\ref{prop:reco}. \qed
\end{proof}

The proof of part~(i) of the following
result mimics that of \cite[Lemma 1.5]{Pl95}.
We shall exploit  the compactness of
the topology~$\T$ of $\overline{{\mathrm{D}}}(\Lalg)$
\cite[Lemma 3.17]{AH93}.

\begin{prop}\label{mpm=ea} Let\, $\Lalg\in \CL$
and let\, $\Y_\T = \overline{{\mathrm{D}}}(\Lalg)$.
\begin{newlist}
\item[{\upshape (i)}]
Every\, $\varphi \in \mpm{\Y_\T}{\twoT_\T}$ is of the
form\, $\bar{e}_a$ for some $a \in L$.
\item[{\upshape (ii)}]
$\Lalg$
is order-isomorphic to\, $\mpm{\Y_\T}{\twoT_\T}$ via the map
$a \mapsto \bar{e}_a$.
\end{newlist}
\end{prop}

\begin{proof}
Since $\varphi$ is $E$-preserving, for any $f \in \varphi^{-1}(1)$ and
$g \in \varphi^{-1}(0)$ we get that $(f,g)\notin E$. Thus there must exist
$a_{fg} \in L$ such that $f(a_{fg})=1$ and $g(a_{fg})=0$. We form the set
of SPH's
$U_{fg} = \{\, h \in \sph{\Lalg}{\twoB} \mid
h(a_{fg}) \neq 0
\,\}$, which is open since
$U_{fg} = \sph{\Lalg}{\twoB} {\setminus} V_{a_{fg}}$. The collection
$\{\, U_{fg} \mid f \in \varphi^{-1}(1) \,\}$ is a cover of $\varphi^{-1}(1)$
since $f \in U_{fg}$ for all $f \in \varphi^{-1}(1)$. Now $\varphi^{-1}(1)$ is compact as
it is a closed subset of the compact space $\overline{{\mathrm{D}}}(\Lalg)$. Thus there is a finite
set $\{\, a_i \mid 1 \le i \le n \,\} \subseteq g^{-1}(0)$ such that
$\varphi^{-1}(1) \subseteq
\{\, \sph{\Lalg}{\twoB}{\setminus} V_{a_i} \mid 1 \le i \le n \,\}$.
We can now set
$a_g = \bigvee \{\, a_i \mid 1 \le i \le n \,\}$. Since $g \in
\sph{\Lalg}{\twoB}$,
we know that
$g^{-1}(0)$ is an ideal and so $g(a_g)=0$.
Furthermore, we have that $\varphi^{-1}(1) \subseteq
\sph{\Lalg}{\twoB}{\setminus}V_{a_g}$.
From Proposition \ref{prop:reco} we have that
$\varphi\in\mpe{\Y}{\twoT}$ and so
using Lemma \ref{3.5+3.7new}(iv) we get $\varphi^{-1}(1) \subseteq W_{a_g}$.

Now $\varphi^{-1}(0)$ is covered by open sets of the form $\sph{\Lalg}{\twoB}{\setminus}W_{a_g}$
for $g~\in~\varphi^{-1}(0)$. Since $\varphi^{-1}(0)$ is also compact, we can take a finite set
$\{\, a^j \mid 1 \le j \le m \,\}$ so that
$\varphi^{-1}(0) \subseteq \bigcup\, \{\, \sph{\Lalg}{\twoB}{\setminus}W_{a^j} \mid 1 \le j \le m \,\}$.
We can also see that
$\varphi^{-1}(1)\subseteq\bigcap\, \{\, W_{a^j} \mid 1 \le j
\le m \,\}$. If we set
$a = \bigwedge \{\, a^j \mid 1 \le j \le m \,\}$
then
$\varphi^{-1}(1) \subseteq W_a$ and $\varphi^{-1}(0) \subseteq \sph{\Lalg}{\twoB} {\setminus} W_a$. Now Lemma
\ref{3.5+3.7new}(iii) gives us that $\varphi^{-1}(0) \subseteq V_a$ and by the maximality of $\varphi$ we
now get $\varphi^{-1}(1)=W_a$ and $\varphi^{-1}(0)=V_a$ and so $\varphi=\bar{e}_a$.

From Proposition \ref{eabar} and part~(i)
we have $\mpm{\Y_\T}{\twoT_\T} = \{\,\bar{e}_a \mid a\in L\,\}$.
Now the map $a \mapsto \bar{e}_a$ is an order-isomorphism
from $\Lalg$ onto
$\mpm{\Y_\T}{\twoT_\T}$ since $a \le b$ if and only if $\bar{e}_a^{-1}(1) \subseteq
\bar{e}_b^{-1}(1)$. To see this, let $a \le b$ and $f\in \bar{e}_a^{-1}(1)$ for $f\in \sph{\Lalg}{\twoB}$. Thus $f(a)=1$ and since
$f^{-1}(1)$ is a filter, we obtain  $f(b)=1$, whence $f\in \bar{e}_b^{-1}(1)$ as required. For the converse, let $a \nleqslant b$ in $\Lalg$. Then $({\uparrow} a,{\downarrow} b)$ is a disjoint filter-ideal pair with an associated special partial homomorphism $f: \Lalg \to \twoB$
with  $f^{-1}(1)= {\uparrow} a$ and $ f^{-1}(0)= {\downarrow} b$. Then $f\in \bar{e}_a^{-1}(1)$ but $f\notin
 \bar{e}_b^{-1}(1)$.
\qed
\end{proof}

We are finally set up to build
a canonical extension of a bounded lattice $\Lalg$,
analogous to the one in Theorem~\ref{canext} but
now based on the larger dual space $\overline{{\mathrm{D}}}(\Lalg)$.
 We  let $Y =  \sph{\Lalg}{\twoB}$ and
shall
show that the complete lattice of maximal partial $E$-preserving
maps from the larger graph $\Y=(Y,E_Y)$
into $\twoT$ is order-isomorphic to the canonical extension from
Theorem \ref{canext}.
To accomplish  this we
recall that $\X = (\mph{\Lalg}{\twoB},E_X)$ and we
define a map $\Psi \colon \mpe{\X}{\twoT} \to \mpe{\Y}{\twoT}$ such that for $\varphi \in \mpe{\X}{\twoT}$
 and $f \in Y$,
\[
(\Psi(\varphi))(f) = \begin{cases}
1 & \text{if } \forall\, g \in \varphi^{-1}(0) \big( (f,g)
\notin E_Y \big),\\
0 & \text{if } \forall\, h \big(\forall\,g \in \varphi^{-1}(0) \big( (h,g) \notin E_Y \Rightarrow (h,f) \notin E_Y \big)\big), \\
- & \text{otherwise.}
\tag{$*$}
\end{cases}
\]
The next proposition establishes, {\it inter alia}, that $\Psi $ is well defined.

\begin{prop}\label{lems_new}
Let\,
$\X = (\mph{\Lalg}{\twoB},E_X)$ and let
 the map \,$\Psi$ be defined by $(*)$ above.
For every \,$\varphi \in \mpe{\X}{\twoT}$,
\begin{newlist}
\item[{\upshape (i)}]
$\varphi^{-1}(1) \subseteq (\Psi(\varphi))^{-1}(1)$ and\,
$\varphi^{-1}(0) \subseteq (\Psi(\varphi))^{-1}(0)$;
\item[{\upshape (ii)}]
the map\, $\Psi(\varphi)$ is a maximal $E$-preserving map
from\, $\Y=(\sph{\Lalg}{\twoB},E_Y)$ into $\twoT$,
and consequently, $\Psi(\varphi)\in \mpe{\Y}{\twoT}$.
\end{newlist}
\end{prop}

\begin{proof}
Consider (i). Let $f\in \mph{\Lalg}{\twoB}$
such that
$f\in\varphi^{-1}(1)$.
By Lemma~\ref{1.3-mod}(ii), for all $g\in\varphi^{-1}(0)$,
we have
$(f,g) \notin E_X$  and so $(f,g) \notin E_Y$.
Hence by the definition of $\Psi(\varphi)$,
we get
$f\in (\Psi(\varphi))^{-1}(1)$.
If $f \in \varphi^{-1}(0)$,
let us
consider $h\in \sph{\Lalg}{\twoB}$ such that
if $g \in \varphi^{-1}(0)$, then $(h,g) \notin E_Y$.
Clearly
$(h,f) \notin E_Y$, which, by the definition of $\Psi(\varphi)$,
gives
$f\in (\Psi(\varphi))^{-1}(0)$.

Now consider~(ii).
First we
assume that $(\Psi(\varphi))(h)=1$ and
$(\Psi(\varphi))(f)=0$ for
some
$h,f \in \sph{\Lalg}{\twoB}$.
Then from
the definition of
$(\Psi(\varphi))^{-1}(1)$
we get that  $(h,g) \notin E_Y$
for every  $g \in \varphi^{-1}(0)$.
 Using this,  we
obtain from the definition of
$(\Psi(\varphi))^{-1}(0)$
that $(h,f) \notin E_Y$.
This shows that
$\Psi(\varphi)$ is $E_Y$-preserving.

Suppose there exists an $E$-preserving map
$\phi\colon\Y\to\twoT$ with $\dom( \Psi(\varphi) )\subsetneq
\dom (\phi)$.
Let $f\in \dom (\phi) \setminus \dom (\Psi(\varphi))$ for
$f \in \sph{\Lalg}{\twoB}$.
First,
let us consider the case $\phi(f)~=~1$.
Since $f \notin (\Psi(\varphi))^{-1}(1)$, there exists $g \in \varphi^{-1}(0)$ such that
$(f,g) \in E_Y$.
By part~(i) we have that $\varphi^{-1}(0) \subseteq (\Psi(\varphi))^{-1}(0)$,
and since $(\Psi(\varphi))^{-1}(0) \subseteq \phi^{-1}(0)$,
we get $g \in \phi^{-1}(0)$. This yields that $\phi$ is not $E_Y$-preserving, a contradiction.
 Now let $\phi(f)=0$. Since $f \notin (\Psi(\varphi))^{-1}(0)$, there exists $g \in (\Psi(\varphi))^{-1}(1)$ with $(g,f) \in E_Y$. But then $g \in \phi^{-1}(1)$ and so $\phi$ is not $E_Y$-preserving, a contradiction. We have shown the maximality of $\Psi(\varphi)$. \qed
\end{proof}

The following lemma will be  deployed in the proof of
Theorem~\ref{XisoY}.

\begin{lem}\label{lem_new3}
Let \,$\Psi $ be defined as in $(*)$.
For\, $\eta \in \mpe {\Y}{\twoT}$,
let\, $\varphi_\eta \colon \X \to \twoT$ be the restriction of\, $\eta$ to\, $\X$.
Then\, $\varphi_\eta \in \mpe{\X}{\twoT}$ and\, $\Psi(\varphi_\eta)=\eta$.
\end{lem}

\begin{proof}
It is easy
to see from the definition that $\varphi_\eta$ is $E_X$-preserving.
Suppose there
exists $\phi \colon \X \to \twoT$ such that
$\dom (\varphi_\eta) \subsetneq \dom (\phi)$.
Consider first the case
where
$f \in \mph{\Lalg}{\twoB}$ with
$f \in \phi^{-1}(1) {\setminus} \varphi_\eta^{-1}(1)$.
This implies that there exists $g \in \eta^{-1}(0)$
such that
$(f,g) \in E_Y$. Consider $h \in \sph{\Lalg}{\twoB}$ defined
by $h^{-1}(1)=f^{-1}(1)$ and $h^{-1}(0) = g^{-1}(0)$. We
can extend $h$ to $h_* \in \mph{\Lalg}{\twoB}$ and by
Lemma~\ref{3.5+3.7new}~(i) we get $h_* \in \varphi_\eta^{-1}(0)$. But
$\varphi_\eta^{-1}(0) \subseteq \phi^{-1}(0)$ and $(f,h) \in E_X$,
contradicting $\phi$ being $E_X$-preserving. Using the same method
we can show that the existence of
$f \in \phi^{-1}(0) {\setminus} \varphi_\eta^{-1}(0)$ implies that
$\phi$ cannot be $E_X$-preserving.
Hence $\varphi_\eta \in \mph{\Lalg}{\twoB}$.

By the definition of the
lattice order of
$\fC{\Y}=\mpe{\Y}{\twoT}$ (see Theorem~\ref{compl}),
it suffices to prove
$(\Psi(\varphi_\eta))^{-1}(1)=\eta^{-1}(1)$.
If $f \in \eta^{-1}(1)$
then for all $g \in \eta^{-1}(0)$, we have $(f,g) \notin E_Y$ where
$f,g\in \sph{\Lalg}{\twoB}$.
In particular this applies
to all $g \in \eta^{-1}(0) \cap \mph{\Lalg}{\twoT} = \varphi_\eta^{-1}(0)$.
Hence $f \in (\Psi(\varphi_\eta))^{-1}(1)$ and so
$\eta^{-1}(1)\subseteq (\Psi(\varphi_\eta))^{-1}(1)$.

If $f \notin \eta^{-1}(1)$, there exists $g \in \eta^{-1}(0)$ such that
$(f,g) \in E_Y$. Let
$h\in \sph{\Lalg}{\twoT}$ be defined by
letting
$h^{-1}(1)=f^{-1}(1)$ and $h^{-1}(0)=g^{-1}(0)$.
We know that we can extend $h$ to
an $h_* \in \mph{\Lalg}{\twoT}$.
Now $g^{-1}(0) \subseteq h_*^{-1}(0)$ and so by Lemma \ref{3.5+3.7new}(i)
we have
that
$h_*~\in~\eta^{-1}(0) \cap \mph{\Lalg}{\twoB} = \varphi_\eta^{-1}(0)$.
Since  $f^{-1}(1) \subseteq h_*^{-1}(1)$,
we have $(f,h_*) \in E_Y$ and so clearly $f \notin (\Psi( \varphi_\eta))^{-1}(1)$.
Thus
$(\Psi(\varphi_\eta))^{-1}(1) \subseteq \eta^{-1}(1)$. \qed
\end{proof}

The next result shows that the complete lattices
$\mpe{\X}{\twoT}$ and $\mpe{\Y}{\twoT}$ are order-isomorphic under
$\Psi$,  and that this map fixes the underlying lattice $\Lalg$.

\begin{thm} \label{XisoY} Let\, $\Lalg\in \CL$.
Let\, $\X = (\mph{\Lalg}{\twoB},E_X)$, $\Y=(\sph{\Lalg}{\twoB},E_Y)$ and
let\, $\Psi$ be the map defined by $(*)$. Then\, $\mpe{\X}{\twoT}$ is
order-isomorphic to\, $\mpe{\Y}{\twoT}$ under the map\, $\Psi$.

Moreover, $\Psi(e_a)=\bar{e}_a$ for any $a \in L$, where
$e_a \in \mpe{\X}{\twoT}$ and\, $\bar{e}_a\in \mpe{\Y}{\twoT}$
are the evaluation maps representing the elements of\,
$\Lalg$ inside its completions\, $\mpe{\X}{\twoT}$ and\, $\mpe{\Y}{\twoT}$, respectively.
\end{thm}

\begin{proof} The surjectivity of $\Psi$
has been proved.
That $\Psi$ is an order-isomorphism can be seen from the fact that $\varphi \le \psi$ if and only if $\varphi^{-1}(1) \subseteq \psi^{-1}(1)$ (see Theorem~\ref{compl}), and this occurs if and only if $\psi^{-1}(0) \subseteq \varphi^{-1}(0)$
(by Lemma~\ref{1.3-mod}),
which, in turn,
occurs if and only if $(\Psi(\varphi))^{-1}(1) \subseteq (\Psi(\psi))^{-1}(1)$ (by the definition of $\Psi$).

The equality $\Psi(e_a)=\bar{e}_a$ follows from
 Lemma~\ref{lem_new3} using the fact that
the restriction of $\bar{e}_a$ to $\X$ is $e_a$. \qed
\end{proof}

Combining the preceding theorem with  Theorem~\ref{canext}
we obtain the following.

\begin{cor}
Let\, $\Lalg \in \CL$.  Then, with respect to the embedding set up by\,
$a \mapsto \bar{e}_a$, the complete lattice\,
$\fC{\Y}=\mpe{\Y}{\twoT}$ is the canonical extension of\, $\Lalg$.
\end{cor}

We are almost ready to define the functors needed
in our categorical approach. First, we make some observations
about the relationship between
the categories we defined earlier:
$\CG$,  of  graphs and $E$-preserving maps,
and its topological
analogue
$\CG_\T$, of
graphs with topology
and
continuous $E$-preserving maps.

From Lemma~\ref{2.1new}
 we note that any $\CL$-graph $\X=(X,E)$
is automatically equipped with
the quasi-orders
$\le_1$ and $\le_2$, as previously defined.
Thus, when considering~$\X$ it is not necessary
explicitly to
state that
$\le_1$ and $\le_2$ are part of the structure of~$\X$.
We note too that on $\twoT$ (which is not an $\CL$-graph) we have quasi-orders
$\le_1$ and $\le_2$  defined as follows:
\[ \le_1 \;=\; E\;=\;\le \quad \text{and} \quad \le_2\;=\;E^{-1}\;=\;\geqslant.\]
We shall
use
the symbol
$\twoT$ to denote both the graph
$(\{0,1\},E)$ and the doubly-ordered set $(\{0,1\},\le_1,\le_2)$.

The quasi-order relations $\le_1$ and $\le_2$ were used by Allwein and Hartonas~\cite{AH93}
as well as by
Urquhart~\cite{U78}.
Furthermore,
Urquhart~\cite{U78} defines
two
maps
$\ell \colon \powerset(X) \to \powerset(X)$ and
$r \colon \powerset(X) \to \powerset(X)$,
given by
\[
\ell(A) = \{\, f \in X \mid \forall\: g \in A (f \nleqslant_1 g) \,\} \quad \text{and}
\quad r(A) = \{\, f \in X \mid \forall\: g \in A (f \nleqslant_2 g) \,\},
\]
for $A \subseteq X$.

The maps $\ell$ and $r$ are order-reversing on
$(\powerset(X),\subseteq)$ and form a Galois connection.
A subset $A \subseteq X$ is defined to be
\defn{$\ell$-stable} if $\ell r (A) =A$ and
\defn{$r$-stable} if $r \ell (A)=A$. The following
result is a simple consequence of the theory of
Galois connections.

\begin{lem} \label{lem:stability} Let\, $\X \in \CY$ and let $A \subseteq X$. Then
\begin{newlist}
\item[{\upshape (i)}] $A$ is $\ell$-stable if and only if\, $A=\ell(B)$ for some\, $B \subseteq X$;
\item[{\upshape (ii)}] $A$ is $r$-stable if and only if\, $A=r(B)$ for some\, $B \subseteq X$.
\end{newlist}
\end{lem}

Now we define morphisms
between $\CL$-graphs.
It is obvious that we
would want any such morphism to preserve the two quasi-orders
$\le_1$ and $\le_2$. We also want morphisms
to preserve the way in which these two
quasi-orders interact with one another.
This we achieve by requiring
that inverse
images of stable subsets of the
codomain  be stable sets in the domain of the morphism.
We state this more precisely below.

Let $\X$ and $\Y$ be $\CL$-graphs
and $\alpha \colon \X \to \Y$. Then $\alpha$ is an
\defn{$\CL$-graph morphism} if
\begin{newlist}
\item[{\upshape (i)}] $\alpha$ preserves $\le_1$ and $\le_2$;
\item[{\upshape (ii)}] if $A \subseteq Y$ is an $\ell$-stable set, then $\alpha^{-1}(r_Y(A))=r_X(\alpha^{-1}(A))$;
\item[{\upshape (iii)}] if $A \subseteq Y$ is an $r$-stable set, then $\alpha^{-1}(l_Y(A))=l_X(\alpha^{-1}(A))$.
\end{newlist}
\noindent
We can now give formal definitions of our categories
$\CY_\T$ and $\CY$:   the former consists of
$\CL$-graphs with  continuous $\CL$-graph morphisms; the latter
is obtained by applying~$^\flat$ to~$\CY_\T$.
Lemma~\ref{Lgraph->Epres} serves to tell us
that
category~$\CY$ (respectively $\CY_\T$) is
a full subcategory of $\CG$
(respectively $\CG_\T$).

\begin{lem} \label{Lgraph->Epres}
Let\, $\X=(X,E_Y)$ and\, $\Y=(Y,E_Y)$ be objects in\, $\CY$ and
consider the map
$\alpha \colon \X \to \Y$. If\, $\alpha$ is an $\CL$-graph morphism,
then $\alpha$ is $E$-preserving.
\end{lem}
\begin{proof}
Let $f,g \in X$ and suppose that $(f,g) \in E_X$.
Then by Lemma~\ref{2.1new}(i)
there exists $h \in X$ such that $f \le_1 h$ and $g \le_2 h$. Since
$\alpha$ is $\le_1$-preserving, we have that
$\alpha(f) \le_1 \alpha(h)$, and since $\alpha$ is
$\le_2$-preserving, we have $\alpha(g) \le_2 \alpha(h)$. Thus there exists
$y=\alpha(h) \in Y$ with the required properties, and so
$(\alpha(f),\alpha(g)) \in E_Y$. \qed
\end{proof}

Not for the first time we are faced with the need to reconcile
two notions which might possibly not coincide.
We shall now prove that MPE's from
$\CL$-graphs
into~$\twoT$ are exactly the maximal partial
$\CL$-graph morphisms into~$\twoT$.

\begin{prop} \label{MPE=MPL}
Let\, $\X=(X,E)$ be an $\CL$-graph, and let\, $\varphi \colon \X \to \twoT$
be a partial map. Then the following are equivalent:
\begin{newlist}
\item[{\upshape (1)}] $\varphi$ is a maximal partial $E$-preserving map;
\item[{\upshape (2)}] $\varphi^{-1}(1)=\ell(\varphi^{-1}(0))$ and\, $\varphi^{-1}(0)=r(\varphi^{-1}(1))$;
\item[{\upshape (3)}] $\varphi$ is a maximal partial $\CL$-graph morphism.
\end{newlist}
\end{prop}
\begin{proof}
We start by showing the equivalence of (1) and (2).
If $\varphi \in \mpe{\X}{\twoT}$ and $g \in \varphi^{-1}(0)$ then, by
Lemma \ref{1.3-mod},
for all $f \in \varphi^{-1}(1)$ we have that $(f,g) \notin E$.
Since $g \le_2 f$ implies that $(f,g)\in E$, we have
for all $f \in \varphi^{-1}(1)$ that $g \nleqslant_2 f$
and hence
$g \in r(\varphi^{-1}(1))$. Thus $\varphi^{-1}(0) \subseteq r(\varphi^{-1}(1))$.
If we suppose that $g \notin \varphi^{-1}(0)$, then by Lemma~\ref{1.3-mod}
we have that there exists $f \in \varphi^{-1}(1)$
such that $(f,g) \in E$. Now by Lemma~\ref{2.1new}(i) there exists $h$
such that $f \le_1 h$ and $g \le_2 h$ and Lemma \ref{3.5+3.7new}(ii) gives us that
$h \in \varphi^{-1}(1)$. Thus $g \notin r(\varphi^{-1}(1))$ and
hence $r(\varphi^{-1}(1)) = \varphi^{-1}(0)$. We can similarly
show that $\ell(\varphi^{-1}(0)) = \varphi^{-1}(1)$.

Now assume (2) and let $A:= \varphi^{-1}(1)$.
To show that $\varphi$ preserves $E$, let $f,g$ be in the domain of $\varphi$, and
suppose $(f,g) \in E$ but $\varphi (f)=1$ and $\varphi(g)=0$. From $\varphi(f)=1$ we get
$f \in \ell r(A) = \ell  (\varphi^{-1}(0))$,
and $\varphi(g)=0$ gives $g \in r(\varphi^{-1}(1)) = r(A)$.
Lemma \ref{2.1new}(i) gives us that there exists $h$
such that
$f\leqslant_1 h $ and $g\leqslant_2 h$.
Now $f \in \ell r(A)$
means that for all
$q \in r(A)$ we have
$f \nleqslant_1 q$
and hence $h \notin r(A)$.
This implies that there exists
$k \in A$ such that $h \le_2 k$. The transitivity of $\le_2$ then implies
that $g \notin r(A)$, a contradiction.
In order to see that $\varphi$ is maximal, suppose that
$\dom \varphi \subseteq \dom \psi$ for some $\psi \in \mpe{\X}{\twoT}$.
It is clear that $A \subseteq \psi^{-1}(1)$ and
$r(A) \subseteq \psi^{-1}(0)$. If
$f \in \psi^{-1}(1)$
and
$g \in \psi^{-1}(0)$, then $(f,g) \notin E$ and hence
$f \nleqslant_1 g$ and $g\nleqslant_2 f$. Thus any $f \in \psi^{-1}(1)$
must be in
$\ell (\psi^{-1}(0))$ and similarly
$\psi^{-1}(0) \subseteq r(\psi^{-1}(1))$.
We then get $\psi^{-1}(1) \subseteq \ell(\psi^{-1}(0)) \subseteq \ell r(A) = A$
and hence $A=\psi^{-1}(1)$.
We also have
that
$r(A) = \psi^{-1}(0)$ and so $\psi = \varphi$,
showing that~$\varphi$ is maximal.
Hence (1) holds.

Now assume (3). Note that $\{0\}$ is $r$-stable and $\{1\}$ is $\ell$-stable in $\twoT$
as $r(\{1\})=\{0\}$ and $\ell(\{0\})=\{1\}$. From (ii) and (iii)
in the definition of morphisms of $\CL$-graphs we get
$\varphi^{-1}(0)= \varphi^{-1}(r(\{1\}))=r(\varphi^{-1}(1))$ and
$\varphi^{-1}(1)=\varphi^{-1}(\ell(\{0\}))=\ell(\varphi^{-1}(0))$.

Finally,
assuming (2), we get that for
any
$f,g \in X$, if
$\varphi(f)=1$ and $\varphi(g)=0$, then
$f \in \ell(\varphi^{-1}(0))$ and
$g \in r(\varphi^{-1}(1))$.
By the definition of
$\ell$ and $r$ we see that
$f \nleqslant_1 g$ and $g \nleqslant_2 f$. Thus
$\varphi$ is $\le_1$- and $\le_2$-preserving on
its domain.
Since
the sets
$\{0\}$ and $\{1\}$
are the $r$- and $\ell$-stable subsets of
$\twoT$, we have that~$\varphi$
obeys properties (ii) and (iii)
required of an $\CL$-graph morphism.
We now suppose that
there exists a partial $\CL$-graph morphism $\psi$ from
$\X$ to $\twoT$ such that
$\dom(\varphi) \subseteq \dom(\psi)$. Suppose there exists
$f \in \dom(\psi)$ but $f \notin \dom(\varphi)$, and
suppose
that
$\psi(f)=1$. Since $f \notin \varphi^{-1}(1)=\ell(\varphi^{-1}(0))$,
there exists $g \in \varphi^{-1}(0) \subseteq \psi^{-1}(0)$ such
that $f \le_1 g$. This implies that $\psi$ is not $\le_1$-preserving.
Similarly, if we suppose that $\psi(f)=0$ we see that
$\psi$ will not be $\le_2$-preserving. Hence $\dom(\varphi)$
is maximal.
So we have shown (3). \qed
\end{proof}

We are, at last, ready to set up the functors we require.
We claim that we can  define $\overline{{\mathrm{D}}}\colon
\CL \to \CY_\T$ and $\overline{{\mathrm{G}}}\colon \CY \to \CL^+$
as follows:

\begin{flalign*}
&\text{on objects:} && \overline{{\mathrm{D}}} \colon  \Lalg \mapsto  (\sph{\Lalg}{\twoB}, E, \T), &   \\
&\text{on morphisms:} && \overline{{\mathrm{D}}} \colon  u \mapsto - \circ u & \\
\intertext{and}
&\text{on objects:} && \overline{{\mathrm{G}}}  \colon  \Y \mapsto  \mpe{\Y}{\twoT}, & \\
&\text{on morphisms:} && \overline{{\mathrm{G}}}  \colon  \alpha  \mapsto - \circ \alpha. &
\end{flalign*}

In the action on morphisms
of
$\overline{\mathrm{D}}$ and
$\overline{\mathrm{G}}$ the domain of the image map
is the set on which the composite map is defined.
In legitimising these definitions what is at stake is
well-definedness:
we must ensure that the images of morphisms
under $\overline{\mathrm{D}}$ and
$\overline{\mathrm{G}}$  are again morphisms for the categories
concerned.

\begin{prop} \label{Doflathom=graphmorph}
Let\, $\Lalg, \K \in \CL$ and let\, $u \colon \Lalg \to \K$ be
a lattice homomorphism. Then
the map
\,
$\smash{\overline{\mathrm{D}}}^\flat(u) \colon
\smash{\overline{\mathrm{D}}}^\flat(\K) \to
\smash{\overline{\mathrm{D}}}^\flat(\Lalg)$
is
an $\CL$-graph morphism.
\end{prop}

\begin{proof} Let $f,g \in \smash{\overline{\mathrm{D}}}^\flat(\K)$ with $f \le_1 g$
and $a \in ((\smash{\overline{\mathrm{D}}}^\flat(u))(f))^{-1}(1)=(f\circ u)^{-1}(1)$.
This gives us that $u(a) \in f^{-1}(1)$ and hence $u(a) \in g^{-1}(1)$.
Thus $a \in (\smash{\overline{\mathrm{D}}}^\flat(u))(g)
^{-1}(1)
$ and so
$\smash{\overline{\mathrm{D}}}^\flat(u)$ is $\le_1$-preserving.
Likewise,
 $\smash{\overline{\mathrm{D}}}^\flat(u)$ is $\le_2$-preserving.

Now consider $A \subseteq \smash{\overline{\mathrm{D}}}^\flat(\Lalg)$ such that
$A$ is $\ell$-stable. We want to show that
$$(\smash{\overline{\mathrm{D}}}^\flat(u))^{-1}(r(A))= r((\smash{\overline{\mathrm{D}}}^\flat(u))^{-1}(A)).$$
Let $f \in \smash{\overline{\mathrm{D}}}^\flat(\K)$ such that
$f \notin r((\smash{\overline{\mathrm{D}}}^\flat(u))^{-1}(A))$. This implies that
there exists $g \in (\smash{\overline{\mathrm{D}}}^\flat(u))^{-1}(A)$ such that
$f \le_2 g$. We note that $g \in (\smash{\overline{\mathrm{D}}}^\flat(u))^{-1}(A)$
if and only if $g \circ u \in A$. Now since
$\smash{\overline{\mathrm{D}}}^\flat(u)$ is $\le_2$-preserving, we have that
$(\smash{\overline{\mathrm{D}}}^\flat(u))(f) = f \circ u \le_2 g \circ u$ and hence
$f \notin (\smash{\overline{\mathrm{D}}}^\flat(u))^{-1}(r(A))$.
Thus we have
$(\smash{\overline{\mathrm{D}}}^\flat(u))^{-1}(r(A)) \subseteq r((
\smash{\overline{\mathrm{D}}}^\flat(u))^{-1}(A))$.

Now if $f \notin (\smash{\overline{\mathrm{D}}}^\flat(u))^{-1}(r(A))$ then
$f \circ u \notin r(A)$ and
hence
there exists $h \in A$ such that $f \circ u \le_2 h$.
Consider the map $g_h \colon \K \to \twoB$
defined for $b \in \K$ by
\[
g_h(b) = \begin{cases}
0 & \text{if } b \in f^{-1}(0), \\
1 & \text{if } b \in {\uparrow}(u[h^{-1}(1)]), \\
- & \text{otherwise.}
\end{cases}
\]
We show that $g_h \in \sph{\Lalg}{\twoB}$.
If $b \in {\uparrow}(u[h^{-1}(1)])$ then there
exists $a \in h^{-1}(1)$ such that $u(a) \le b$.
Since $f \circ u \le_2 h$ we have that
$h^{-1}(1) \cap (f \circ u)^{-1}(0)=\emptyset$ and
so $a \notin (f \circ u)^{-1}(0)$. If $b \in f^{-1}(0)$
then $u(a) \in f^{-1}(0)$, a contradiction.
Hence $g_h^{-1}(1) \cap g_h^{-1}(0) = \emptyset$.
We need to show that $g_h^{-1}(1)$ is a filter of~$\K$,
so we consider $b_1,b_2 \in g_h^{-1}(1)$.
By the definition of $g_h$,
there must exist $a_1,a_2 \in h^{-1}(1)$ such that $u(a_1) \le b_1$
and $u(a_2) \le b_2$. Now since $h^{-1}(1)$ is a filter of~$\Lalg$,
we have that $a_1 \wedge a_2 \in h^{-1}(1)$ and
since $u(a_1 \wedge a_2)$ is a lower bound for $\{b_1, b_2\}$ we have
that $u(a_1 \wedge a_2) \le b_1 \wedge b_2$
and so $b_1 \wedge b_2 \in g_h^{-1}(1)$.
Hence $g_h \in \sph{\K}{\twoB}$.

If $a \in h^{-1}(1)$ then $u(a) \in u[h^{-1}(1)]$ and so $u(a) \in g_h^{-1}(1)$.
Hence $g_h(u(a))=1$ and we have $h \le_1 g_h \circ u$. Now
since $h \in A = \ell r (A)$ and using the definition of $\ell$
and the transitivity of $\le_1$,
we have
that $g_h \circ u \in \ell r(A)=A$.
Now clearly $f \le_2 g_h$ and since
$g_h \in (\smash{\overline{\mathrm{D}}}^\flat(u))^{-1}(A)$
we have $f \notin r((\smash{\overline{\mathrm{D}}}^\flat(u))^{-1}(A))$.
Thus $r((\smash{\overline{\mathrm{D}}}^\flat(u))^{-1}(A)) \subseteq (\smash{\overline{\mathrm{D}}}^\flat(u))^{-1}(r(A))$. \qed
\end{proof}

In order to show the well-definedness of $\overline{\mathrm{G}}$,
we must verify that the image of a morphism under $\overline{\mathrm{G}}$ takes MPE's to MPE's.
To address this  question we first prove a technical lemma.

\begin{lem} \label{lem:gmaps}
Let\, $\X,\Y \in \CY$ and let\, $\alpha \colon \X \to \Y$ be such
that\, $\alpha \in \CY(\X,\Y)$. Further, let\,
$\varphi~\in~\mpe{\Y}{\twoT}$
and $f \in X$. Then
\begin{newlist}
\item[{\upshape (i)}] if there exists $m \in \varphi^{-1}(0)$ such that\, $(\alpha(f),m) \in E_Y$
then there exists $g_m \in X$ such that\, $\varphi(\alpha(g_m))=0$ and\, $(f,g_m) \in E_X$;
\item[{\upshape (ii)}] if there exists $n \in \varphi^{-1}(1)$ such that\, $(n,\alpha(f)) \in E_Y$ then
there exists $g_n \in X$ such that\, $\varphi(\alpha(g_n))=1$ and\, $(g_n,f) \in E_X$.
\end{newlist}
\end{lem}
\begin{proof} Considering (i), we have from Lemma~\ref{1.3-mod} that
$\alpha(f)
\notin \varphi^{-1}(1)$.
This implies that
$f
\notin \alpha^{-1}(\varphi^{-1}(1))=\alpha^{-1}(\ell(\varphi^{-1}(0)))$. Now since
$\alpha^{-1}$ preserves $\ell$-stable sets,
we have that
$f
\notin \ell (\alpha^{-1}(\varphi^{-1}(0)))$. This implies that
there exists $g_m \in \alpha^{-1}(\varphi^{-1}(0))$ such that $f \le_1 g_m$. Now clearly
$\alpha(g_m)\in \varphi^{-1}(0)$ and $(f,g_m) \in E_X$. \qed
\end{proof}

\begin{prop} \label{Gmax}
Let\, $\X,\Y \in \CY$ and let\, $\alpha \colon \X \to \Y$ be
an $\CL$-graph morphism.
Then for any\, $\varphi \in
\mpe{\Y}{\twoT}$
it is the case that\,
$(\overline{\mathrm{G}}(\alpha))(\varphi) \in \mpe{\X}{\twoT}$.
\end{prop}
\begin{proof}
We have from Lemma \ref{Lgraph->Epres} that $\alpha$ is
$E$-preserving, and hence we can apply
Lemma~\ref{E-pres}
to conclude that
$(\overline{\mathrm{G}}(\alpha))(\varphi)$ is a
partial
$E$-preserving map from $\X$ to $\twoT$.
We show that its domain is maximal.

Suppose that there exists an $E$-preserving map $\psi$
extending
$(\overline{\mathrm{G}}(\alpha))(\varphi)$ and
then suppose
that
$f \notin \dom
((\overline{\mathrm{G}}(\alpha))(\varphi))$
for some
$f \in X$.
We want to show that
$f \notin \dom (\psi)$.
Since $f \notin \dom
(\varphi \circ \alpha)$
we know that
$(\varphi \circ \alpha)(f) \neq 1$
and hence by Lemma~\ref{1.3-mod}
there must exist $m_0 \in \varphi^{-1}(0)$ such that
$(\alpha(f),m_0) \in E_Y$.
Similarly, since
$f \notin (\varphi \circ \alpha)^{-1}(0)$
we know that
$\varphi (\alpha(f)) \neq 0$
and hence by Lemma~\ref{1.3-mod}
there
exists
$m_1 \in \varphi^{-1}(1)$ such that
$(m_1,\alpha(f)) \in E_Y$.

Now from Lemma \ref{lem:gmaps} we
get
$g_{m_0}, g_{m_1}\in
\X$
such that
$\varphi(\alpha(g_{m_1}))=1$ with
$(g_{m_1},f) \in E_X$ and
$\varphi(\alpha(g_{m_0}))=0$
with $(f,g_{m_0}) \in E_X$.
From this we see that $g_{m_1} \in
((\overline{\mathrm{G}}(\alpha))(\varphi))^{-1}(1)$
whence $\psi(g_{m_1})=1$ and $g_{m_0} \in
((\overline{\mathrm{G}}(\alpha))(\varphi))^{-1}(0)$
whence $\psi(g_{m_0})=0$.
If we now suppose, for contradiction, that $f \in \dom (\psi)$, then
from $(g_{m_1},f) \in E_X$ we obtain $1=\psi(g_{m_1}) \le \psi(f)$ and from
$(f,g_{m_0}) \in E_X$ we obtain $\psi(f) \le \psi(g_{m_0})=0$, which is impossible. \qed
\end{proof}

The next result confirms
that the codomain of $\overline{{\mathrm{G}}}$ really is $\CL^+$, as our commuting diagram
demands.  The theorem which follows shows that under $^\delta$ an $\CL$-morphism
lifts to an
$\CL^+$-morphism.
We emphasise  that the  proof of the first result
is a completely routine
definition chase.
By contrast,
lattice homomorphisms were regarded  by Gehrke and Harding \cite{GH01}
as instances
of additional operations which are
join- and meet-preserving, and their liftings treated using the machinery
developed to handle extensions of maps in general.

\begin{prop} \label{G-on-morphs} Let \,$\X, \Y \in \CY$ and $\alpha \in \CY(\X,\Y)$.
Then\, $\overline{\mathrm{G}}(\alpha) \colon \mpe{\Y}{\twoT} \to \mpe{\X}{\twoT}$
is a complete lattice homomorphism.
That is,
$\overline{\mathrm{G}}(\alpha) \in \CL^+(\overline{\mathrm{G}}(\Y),\overline{\mathrm{G}}(\X))$.
\end{prop}

\proof We prove that $\overline{\mathrm{G}}(\alpha)$ is meet-preserving. The fact
that $\overline{\mathrm{G}}(\alpha)$ is join-preserving will
follow by a similar argument.
Let $f \in X$ and consider the collection
$\{\, \varphi_i \mid i \in I\,\} \subseteq \mpe{\Y}{\twoT}$. Then
\begin{align*}
\allowdisplaybreaks
\qquad \quad \;\;\, f \in \big( \DM_{i \in I}
(\overline{\mathrm{G}} (\alpha)(\varphi_i)) \big)^{-1}(1) & \:\Longleftrightarrow\: f \in \bigcap_{i \in I}
(\overline{\mathrm{G}}(\alpha)(\varphi_i))^{-1}(1) \\
& \:\Longleftrightarrow\: (\forall\, i \in I)\; (\varphi_i \circ \alpha)(f) = 1 \\
& \:\Longleftrightarrow\: (\forall\, i \in I)\; \varphi_i (\alpha(f)) = 1 \\
& \:\Longleftrightarrow\: \alpha(f) \in \bigcap_{i \in I} \varphi_i^{-1}(1) \\
& \:\Longleftrightarrow\: \big( \DM_{i \in I} \varphi_i \big)(\alpha(f))=1 \\
& \:\Longleftrightarrow\: \big( (\DM_{i \in I} \varphi_i) \circ \alpha \big)(f)=1 \\
& \:\Longleftrightarrow\: f \in \big( \overline{\mathrm{G}}(\alpha) (\DM_{i \in I} \varphi_i ) \big)^{-1}(1). \quad \quad & \end{align*}
\qed

We can now present the main result of this section.

\begin{thm}
\label{hom}
Let $u \in \CL(\Lalg, \K)$ be a homomorphism of bounded lattices. Then
\[ \overline{\mathrm{G}} \hspace*{.04cm} \smash{\overline{{\mathrm{D}}}}^\flat(u) \colon
\overline{\mathrm{G}}\hspace*{.04cm}\smash{\overline{{\mathrm{D}}}}^\flat (\Lalg) \to \overline{\mathrm{G}}
\hspace*{.04cm}\smash{\overline{{\mathrm{D}}}}^\flat (\K) \]
given by the
composition of functors\,
$\overline{{\mathrm{G}}} \circ ^\flat
\circ
 {\overline{\mathrm{D}}} $ is a
complete homomorphism of the corresponding canonical extensions.
\end{thm}
\begin{proof} This follows from Propositions~\ref{Doflathom=graphmorph} and~\ref{G-on-morphs}. \qed
\end{proof}

In conclusion, we
sum up what we have achieved.
\begin{thm} Let the categories $\CL$, $\CL^+$,
$\CY_\T$
and $\CY$
be as  above and construct\, $\overline{{\mathrm{D}}}$ and\, $\overline{{\mathrm{G}}}$
as indicated.  Then
\begin{newlist}
\item[{\upshape (i)}]  $\overline{{\mathrm{D}}}\colon \CL \to \CY_\T$ and\,
$\overline{{\mathrm{G}}}\colon \CY \to \CL^+$ are well-defined functors;
\item[{\upshape (ii)}] the functor\, $\overline{{\mathrm{D}}}$ has a right adjoint\,
$\overline{\mathrm{E}}$ such that the unit of the adjunction
is given by evaluation maps which are isomorphisms;

\item[{\upshape (iii)}] with the categories and functors as defined above,
the diagram in Fig.~{\upshape \ref{fig:Lsquare}} commutes, so that the canonical extension
functor on~$\CL$ factorises as the composition\, $\overline{{\mathrm{G}}} \circ ^\flat \circ
\overline{{\mathrm{D}}}$.
\end{newlist}
\end{thm}

\begin{acknowledgements}
The authors express their thanks to the referee as well as to Miroslav Plo\v s\v cica for their useful comments with respect to the proof of Proposition~\ref{compact}.
\end{acknowledgements}

\bibliographystyle{spmpsci}

\end{document}